%% file: cliffordx.tex
\documentclass[12pt]{book}
\usepackage{here,amsfonts,amssymb,makeidx}
\usepackage{amsmath,amsthm}
\def\reals{\mathbb{R}}
\def\complex{\mathbb{C}}
\def\integs{\mathbb{Z}}
\def\natnums{\mathbb{N}}

\def\mapdef#1#2#3{#1\co #2\rightarrow #3}
\def\Ker{\mathrm{Ker}\,}
\def\dimm{\mathrm{dim}}

\def\Hom#1#2#3{\mathrm{Hom}_{#1}(#2, #3)}

\def\buildrell#1\over#2{\mathrel{\mathop{\null#1}\limits_{#2}}}

\def\minibull{\bullet}

\def\Pin{{\bf Pin}} 
\def\Spin{{\bf Spin}} 
\def\mPin#1{{\bf Pin}(#1)} 
\def\mSpin#1{{\bf Spin}(#1)} 
\def\mCl{\mathrm{Cl}}
\def\mH#1{\mathbf{H}(#1)}
\def\mU#1{\mathbf{U}(#1)}
\def\mSU#1{\mathbf{SU}(#1)}
\def\mSO#1{\mathbf{SO}(#1)}
\def\mO#1{\mathbf{O}(#1)}
\def\quat{\mathbb{H}}
\def\mGL#1{\mathbf{GL}(#1)}
\def\qone{\mathbf{1}}
\def\qi{\mathbf{i}}
\def\qj{\mathbf{j}}
\def\qk{\mathbf{k}}
\def\mfrac#1{{\mathfrak{#1}}}
\def\id{\mathrm{id}}
\def\novect#1{#1}
\def\norme#1{\left\|#1\right\|}
\def\transpos#1{#1^{\top}}
\def\amsmata#1#2#3#4{%
\begin{pmatrix}
#1 & #2\\
#3 & #4
\end{pmatrix}
}

\input xy
\xyoption{all}

\newbox\sample
\input lamacb2-ams.tex

\input mac-new.tex

\title{Clifford Algebras, Clifford Groups,\\
and a Generalization of the Quaternions:\\
The \Pin\ and \Spin\ Groups
}
\author{Jean Gallier\\
Department of Computer and Information Science\\
University of Pennsylvania\\
Philadelphia, PA 19104, USA\\
e-mail: {\tt jean@cis.upenn.edu}\\
\ \\
}
\begin{document}
\maketitle
\input abstcliff.tex
\tableofcontents
\vfill\eject
\input clifford1.tex

\input cliffordx.bbl
\bibliographystyle{plain} 
\end{document}

%% file: lamacb2-ams.tex
\theoremstyle{plain}
\newtheorem{thm}{Theorem}[chapter]
\newtheorem{theorem}[thm]{Theorem}
\newtheorem{lemma}[thm]{Lemma}
\newtheorem{cor}[thm]{Corollary}

\newtheorem{prop}[thm]{Proposition}
\newtheorem{proposition}[thm]{Proposition}

\theoremstyle{definition}
\newtheorem{example}{Example}[chapter]

\newtheorem{definition}{Definition}[chapter]
\newtheorem{defin}[definition]{Definition}

%% file: mac-new.tex
\newbox\sample
\newif\ifproofmode
\newif\ifsymindex
\global\symindexfalse
\newwrite\inx
\def\indsyma#1#2{\ifproofmode\marginpar{$\scriptstyle#1$}\fi%
\ifx#2\empty\write\inx{$\noexpand#1$,\space\thepage}%
\write\inx{\string\newline}\else%
\write\inx{$\noexpand#1$,\space#2,\space\thepage}%
\write\inx{\string\newline}\fi\ignorespaces}%
\def\indsym#1#2{\ifsymindex%
\ifproofmode\marginpar{$\scriptstyle#1$}\fi%
\ifx#2\empty\write\inx{\string\item \space$\noexpand#1$,\space\thepage}%
\else%
\write\inx{\string\item \space$\noexpand#1$,\space#2,\space\thepage}%
\fi\ignorespaces\fi}%
\newskip\dangerskipb
\newskip\dangerskip
\def\hang{\hangindent\dangerskip}

\def\tensor{\otimes}

\def\remark{\bigskip\noindent{\bf Remark:}\enspace}

\def\remarks{\bigskip\noindent{\bf Remarks:}\enspace}

\font\manual=manfnt at 12pt
\def\danbend{{\manual\char127}}
\def\datanger{\medbreak\begingroup\clubpenalty=10000
 \def\par{\endgraf\endgroup\medbreak} \noindent\hang\hangafter=-2
 \hbox to0pt{\hskip-3.5pc\danbend\hfill}}
\outer\def\danger{\datanger}%
\def\ddatanger{\medbreak\begingroup\clubpenalty=10000
 \def\par{\endgraf\endgroup\medbreak} \noindent\hang\hangafter=-2
 \hbox to0pt{\hskip-3.5pc\danbend\kern1pt%
\danbend\hfill}}

\def\dobdownarrow{\mathop{\vbox{\kern2pt \hbox{$\Big\downarrow$}\kern-16.5pt
                          \nointerlineskip\hbox{$\Big\downarrow$}}}}

\def\lrightarrow{\hbox to 25pt{\rightarrowfill}}

\def\supexp{exp(m,n,p)=m^{m^{m^{\cdot^{\cdot^{\cdot^{m^{p}}}}}}}
\vbox{\hbox{$\Big\}\scriptstyle n$}\kern0pt}}

\def\supexpo#1#2#3{#1^{#1^{\cdot^{\cdot^{\cdot^{#1^{#2}}}}}}
\vbox{\hbox{$\Big\}\scriptstyle #3$}\kern0pt}}

\def\sqr#1#2{{\vcenter{\hrule height .#2pt
         \hbox{\vrule width.#2pt height#1pt \kern#1pt
             \vrule width.#2pt}
         \hrule height.#2pt}}}

\def\co{\colon}

\newskip\bogcentering \bogcentering= 0pt plus 1000pt minus 1000pt 

\def\matth{\mathsurround=0pt}
\def\fakrightarrowfill{$\matth \mathord- \mkern-6mu
  \cleaders\hbox{$\mkern-2mu \mathord- \mkern-2mu$}\hfill
 \mkern-6mu \mathord\rightarrow$}

\def\fakoverrightarrow#1{\vbox{\ialign{##\crcr
  \fakrightarrowfill\crcr\noalign{\kern-1pt\nointerlineskip}
 $\hfil\displaystyle{#1}\hfil$\crcr}}}

\def\oldcases#1{\left\{\,\vcenter{\normalbaselines\matth
  \ialign{$##\hfil$&\quad##\hfil\crcr#1\crcr}}\right.}

\newif\ifdtatp

\def\displaty{%
\global \dtatptrue \openup \jot \matth \everycr{\noalign{\ifdtatp \global 
\dtatpfalse \vskip -\lineskiplimit \vskip \normallineskiplimit \else 
\penalty \interdisplaylinepenalty \fi }}}

\def\displaylignes#1{\displaty
   \halign{\hbox to\displaywidth{$\displaystyle##$}\crcr
   #1\crcr}}
\def\leqaligneno#1{\displaty \tabskip=\bogcentering
 \halign to\displaywidth{\hfil$\displaystyle{##}$\tabskip=0pt
 &$\displaystyle{{}##}$\hfil\tabskip=\bogcentering
 &\kern-\displaywidth\rlap{$##$}\tabskip=\displaywidthpt\crcr
 #1\crcr}}
\def\ligne{\hbox to\hsize}
\newdimen\nouvpagewidth
\newdimen\offwidth
\newdimen\lawidthoui
\nouvpagewidth=\lawidthoui
\def\kboxit#1{\vbox{\hrule\hbox{\vrule\kern3pt
              \vbox{\kern3pt#1\kern3pt}\kern3pt\vrule}\hrule}}

\def\kboxitb#1{\vbox{\hrule\hbox{\vrule\kern3pt
              \vbox{\kern3pt#1\kern3pt}\kern3pt\vrule}\hrule}}

\def\laboxaround#1{
\aboxaround{\hbox to\hsize{\hfill\box2\hfill}}{#1}
}

\def\boxar#1#2{
\aboxaround{\hbox to\hsize{\hfill#1\hfill}}{#2}
}

\def\aboxaround#1#2{
\setbox4=\vbox{\hsize #2\noindent\strut#1\strut}
\kboxitb{\box4}}

\def\kframeit#1{\vbox{\hrule\hbox{\vrule\kern5pt
              \vbox{\kern5pt#1\kern5pt}\kern5pt\vrule}\hrule}}

\newskip\savnormalbaselineskip
\newskip\savnormallineskip
\newdimen\savnormallineskiplimit

%% file: abstcliff.tex
\ \vfill\eject
\begin{center}
{\large \bf
Clifford Algebras, Clifford Groups,\\
and a Generalization of the Quaternions:\\
The \Pin\ and \Spin\ Groups
}
\end{center}

\vspace{1cm}
\begin{center}
Jean Gallier
\end{center}

\vspace{2cm}

\noindent
{\bf Abstract:\/}
One of the main goals of these notes is to explain
how rotations in $\reals^n$ are induced by the action
of a certain group $\mSpin{n}$ on $\reals^n$,
in a way that generalizes the action of the unit
complex numbers $\mU{1}$ on $\reals^2$, and
the action of the unit quaternions $\mSU{2}$
on $\reals^3$ ({\it i.e.\/}, the action
is defined in terms of multiplication in a larger
algebra containing both the group $\mSpin{n}$ and $\reals^n$). 
The group $\mSpin{n}$, called a
{\it spinor group\/}, is defined
as a certain subgroup of units of an algebra
$\mCl_n$, the {\it Clifford algebra\/} associated
with $\reals^n$.

\medskip
Since the spinor groups are certain well chosen
subgroups of units of Clifford algebras,
it is necessary to investigate
Clifford algebras to get a firm understanding
of spinor groups. These notes provide a tutorial
on Clifford algebra and the groups $\Spin$ and $\Pin$,
including a study of
the structure of the Clifford algebra
$\mCl_{p, q}$ associated with a nondegenerate symmetric
bilinear form of signature $(p, q)$ and culminating
in the beautiful ``$8$-periodicity theorem'' of
Elie Cartan and Raoul Bott (with proofs).

%% file: clifford1.tex
\chapter[Clifford Algebras, Clifford Groups, $\Pin$  and $\Spin$]
{Clifford Algebras, Clifford Groups, and the Groups
$\mPin{n}$ and $\mSpin{n}$}
\label{chapcliff1}
\section{Introduction: Rotations As Group Actions}
\label{seccliff0}
The main goal of this chapter is to explain
how rotations in $\reals^n$ are induced by the action
of a certain group $\mSpin{n}$ on $\reals^n$,
in a way that generalizes the action of the unit
complex numbers $\mU{1}$ on $\reals^2$, and
the action of the unit quaternions $\mSU{2}$
on $\reals^3$ ({\it i.e.\/}, the action
is defined in terms of multiplication in a larger
algebra containing both the group $\mSpin{n}$ and $\reals^n$). 
The group $\mSpin{n}$, called a
{\it spinor group\/}, is defined
as a certain subgroup of units of an algebra
$\mCl_n$, the {\it Clifford algebra\/} associated
with $\reals^n$. Furthermore, for $n\geq 3$,
we are lucky, because
the group $\mSpin{n}$ is topologically
simpler than the group $\mSO{n}$. Indeed, for $n\geq 3$,
the group  $\mSpin{n}$ is simply connected
(a fact that it not so easy to prove without
some machinery),
whereas $\mSO{n}$ is not simply connected.
Intuitively speaking, $\mSO{n}$ is more twisted
than $\mSpin{n}$. In fact, we will see that
$\mSpin{n}$ is a double cover of $\mSO{n}$.

\medskip
Since the spinor groups are certain well chosen
subroups of units of Clifford algebras,
it is necessary to investigate
Clifford algebras to get a firm understanding
of spinor groups. This chapter provides a tutorial
on Clifford algebra and the groups $\Spin$ and $\Pin$,
including a study of
the structure of the Clifford algebra
$\mCl_{p, q}$ associated with a nondegenerate symmetric
bilinear form of signature $(p, q)$ and culminating
in the beautiful ``$8$-periodicity theorem'' of
Elie Cartan and Raoul Bott (with proofs).
We also explain when $\mSpin{p,q}$ is a double-cover
of $\mSO{p,q}$. The reader should be warned that 
a certain amount of algebraic (and topological) background
is expected.
This being said, perseverant readers will
be rewarded by being exposed to some beautiful 
and nontrivial concepts and results, 
including Elie Cartan and Raoul Bott 
``$8$-periodicity theorem.'' 

\medskip
Going back to rotations as transformations induced
by group actions, recall that if 
$V$ is a vector space, a {\it linear action (on the left) of
a group $G$ on $V$\/} is a map $\mapdef{\alpha}{G\times V}{V}$
satisfying the following conditions,
where, for simplicity of notation, we denote $\alpha(g, v)$
by $g\cdot v$:
\begin{enumerate}
\item[(1)]
$g\cdot (h \cdot v) = (gh)\cdot v$, for all $g, h\in G$ and $v\in V$;
\item[(2)]
$1\cdot v = v$, for all $v\in V$, where $1$ is the identity of the group $G$;
\item[(3)]
The map $v \mapsto g\cdot v$ is a linear isomorphism of $V$ for
every $g\in G$.
\end{enumerate}

For example, the (multiplicative) group $\mU{1}$
of unit complex numbers acts on $\reals^2$ (by identifying
$\reals^2$ and $\complex$) {\it via\/} complex multiplication:
For every $z = a + ib$ (with $a^2 + b^2 = 1$),
for every $(x, y) \in \reals^2$ (viewing $(x, y)$ as
the complex number $x + iy$), 
\[z\cdot (x, y) = (ax - by, ay + bx).\]
Now, every unit complex number is of the form
$\cos\theta + i\sin\theta$, and thus
the above action of $z = \cos\theta + i\sin\theta$ on
$\reals^2$ corresponds to the rotation of angle $\theta$
around the origin.
In the case $n = 2$, the groups $\mU{1}$ and $\mSO{2}$ are isomorphic,
but this is an exception.

\medskip
To represent rotations in $\reals^3$ and $\reals^4$, we need the
quaternions. For our purposes, it is convenient to define the
quaternions as certain $2\times 2$ complex matrices.
Let $\mathbf{1}, \mathbf{i}, \mathbf{j}, \mathbf{k}$ be the 
    matrices
\begin{align*}
\mathbf{1} & = 
\begin{pmatrix}
1 & 0 \\
0 & 1
\end{pmatrix},
&  
\mathbf{i} & = 
\begin{pmatrix}
i & 0 \\
0 & -i
\end{pmatrix}, 
&
\mathbf{j} & = 
\begin{pmatrix}
0 & 1 \\
-1 & 0
\end{pmatrix},
&  
\mathbf{k} & = 
\begin{pmatrix}
0 & i \\
i & 0
\end{pmatrix}, 
\end{align*}
and let $\mathbb{H}$ be the set of all matrices of the form
\[
X = a \mathbf{1} + b \mathbf{i}  + c \mathbf{j} + d \mathbf{k}, \quad
a, b, c, d \in \reals.  
\]
Thus, every matrix in $\mathbb{H}$ is of the form
\[
X =
\begin{pmatrix}
a + i b & c + i d \\
-(c - i d)  & a - i b
\end{pmatrix}, \quad a, b, c, d\in \reals.
\]
The quaternions  $\mathbf{1}, \mathbf{i}, \mathbf{j}, \mathbf{k}$
satisfy the famous identities discovered by Hamilton:
\begin{align*}
& \mathbf{i}^2 =  \mathbf{j}^2 = \mathbf{k}^2 =
\mathbf{i}\mathbf{j}\mathbf{k} = -\mathbf{1}, \\
& \mathbf{i}\mathbf{j}    = -\mathbf{j}\mathbf{i} = \mathbf{k}, \\
& \mathbf{j}\mathbf{k}    = -\mathbf{k}\mathbf{j} = \mathbf{i}, \\
& \mathbf{k}\mathbf{i}    = -\mathbf{i}\mathbf{k} = \mathbf{j}.
\end{align*}
As a consequence, it can be verified that
$\mathbb{H}$ is a skew field (a noncommutative field) called
the {\it quaternions\/}.  It is also a real vector space of dimension
$4$ with basis  $(\mathbf{1}, \mathbf{i}, \mathbf{j}, \mathbf{k})$;
thus as a vector space, $\mathbb{H}$ is isomorphic to $\reals^4$.
The {\it unit quaternions\/} are the
quaternions such that
\[
\det(X) = a^2 + b^2 + c^2 + d^2 = 1.
\]
Given any quaternion
$X = a \mathbf{1} + b \mathbf{i}  + c \mathbf{j} + d \mathbf{k}$,
the {\it conjugate\/} $\overline{X}$ of $X$ is given by
\[
\overline{X} = a \mathbf{1} - b \mathbf{i}  - c \mathbf{j} - d \mathbf{k}.
\] 
It is easy to check that the matrices associated with the unit
quaternions are exactly the matrices in $\mSU{2}$. Thus, we call 
$\mSU{2}$ the group of unit quaternions.

\medskip
Now we can  define an action of the group of
unit quaternions $\mSU{2}$ on $\reals^3$.
For this, we use the fact that $\reals^3$ can
be identified with the pure quaternions in $\quat$,
namely, the quaternions of the form
$x_1\qi + x_2\qj + x_3\qk$, where $(x_1, x_2, x_3)\in \reals^3$.
Then, we define the action of $\mSU{2}$ over
$\reals^3$ by
\[Z \cdot X =  ZXZ^{-1} = ZX\overline{Z},\]
where $Z\in \mSU{2}$ and $X$ is any pure quaternion.
Now, it turns out that the map $\rho_Z$
(where $\rho_Z(X) = ZX\overline{Z}$)
is indeed a rotation, and that the map
$\rho\co Z \mapsto \rho_Z$ is 
a surjective homomorphism
$\mapdef{\rho}{\mSU{2}}{\mSO{3}}$ whose kernel is
$\{-\qone, \qone\}$, where 
$\qone$ denotes the multiplicative unit quaternion.
(For details, see Gallier \cite{Gallbook2}, Chapter 8).

\medskip
We can also define an action of the group
$\mSU{2}\times \mSU{2}$ over $\reals^4$, by identifying
$\reals^4$ with the quaternions. In this case,
\[(Y, Z) \cdot X = YX\overline{Z},\]
where $(Y, Z) \in \mSU{2}\times \mSU{2}$ and
$X\in \quat$ is any quaternion.
Then, the map $\rho_{Y, \overline{Z}}$ is a rotation
(where $\rho_{Y, \overline{Z}}(X) = YX\overline{Z}$), and the map
$\rho\co (Y, Z) \mapsto \rho_{Y, \overline{Z}}$ 
is a surjective homomorphism
$\mapdef{\rho}{\mSU{2}\times \mSU{2}}{\mSO{4}}$ 
whose kernel is 
$\{(\qone, \qone), (-\qone, -\qone)\}$.
(For details, see Gallier \cite{Gallbook2}, Chapter 8).

\medskip
Thus, we observe that for $n = 2, 3, 4$, the
rotations in $\mSO{n}$ can be realized {\it via\/}
the linear action of some group
(the case $n = 1$ is trivial, since $\mSO{1} = \{1, -1\}$).
It is also the case that the action of
each group can be somehow be described in terms of multiplication
in some larger algebra ``containing'' the original vector
space $\reals^n$
($\complex$ for $n = 2$, $\quat$ for $n = 3, 4$).
However, these
groups appear to have been discovered in an ad hoc
fashion, and there does not appear to be any universal
way to define the action of these groups on $\reals^n$.
It would certainly be nice if the action was always of the form
\[Z \cdot X =  ZXZ^{-1} (= ZX\overline{Z}).\]
A systematic way of constructing groups 
realizing  rotations in terms of linear action,
using a uniform notion of action, does exist.
Such groups are the spinor groups,
to be described in the following sections.

\section{Clifford Algebras}
\label{seccliff1}
We explained in Section \ref{seccliff0}  how
the rotations in $\mSO{3}$ can be realized
by the linear action of the group of unit quaternions
$\mSU{2}$ on $\reals^3$, and how the rotations
in  $\mSO{4}$ can be realized 
by the linear action of the group
$\mSU{2}\times \mSU{2}$
on $\reals^4$.

\medskip
The main reasons why the rotations in $\mSO{3}$
can be represented by unit quaternions are the following:
\begin{enumerate} 
\item[(1)]
For every nonzero vector $u\in\reals^3$,
the reflection $s_u$ about the hyperplane
perpendicular to $u$ is represented by the map
\[v \mapsto -u v u^{-1},\]
where $u$ and $v$ are viewed as pure quaternions in $\quat$
({\it i.e.}, if $u = (u_1, u_2, u_2)$, then
view $u$ as $u_1\qi + u_2\qj + u_3\qk$, and similarly for $v$). 
\item[(2)]
The group $\mSO{3}$ is generated by the reflections.
\end{enumerate} 

As one can imagine, a successful generalization
of the quaternions, {\it i.e.}, the discovery of a group
$G$ inducing the rotations in $\mSO{n}$ {\it via\/}
a linear action, 
depends on the ability to generalize properties
(1) and (2) above. Fortunately, it is true that
the group $\mSO{n}$ is generated by the hyperplane
reflections. In fact, this is also true for the 
orthogonal group $\mO{n}$, and more generally
for the group of direct isometries $\mO{\Phi}$ of any nondegenerate
quadratic form $\Phi$, by the {\it Cartan-Dieudonn\'e theorem\/}
(for instance, see Bourbaki \cite{BourbakiA3}, or 
Gallier \cite{Gallbook2}, Chapter 7, Theorem 7.2.1).
In order to generalize (2), we need to understand
how the group $G$ acts on $\reals^n$.
Now, the case $n = 3$ is special, because
the underlying space $\reals^3$ on which 
the rotations act can be embedded as the pure quaternions
in $\quat$. The case $n = 4$ is also special,
because $\reals^4$ is the underlying space of $\quat$.
The generalization
to $n\geq 5$ requires more machinery, namely,
the notions of Clifford groups and Clifford algebras.

\medskip
As we will see, for every $n \geq 2$, there is
a compact, connected (and simply connected when $n\geq 3$) group
$\mSpin{n}$, the ``spinor group,''
and a surjective homomorphism
$\mapdef{\rho}{\mSpin{n}}{\mSO{n}}$ whose kernel
is $\{-1, 1\}$. 
This time, $\mSpin{n}$ acts directly on $\reals^n$,
because $\mSpin{n}$ is a certain subgroup of
the group of units of the {\it Clifford algebra\/}
$\mCl_n$, and $\reals^n$ is naturally a subspace
of $\mCl_n$. 

\medskip
The group of unit quaternions $\mSU{2}$ turns out to be isomorphic to
the spinor group  $\mSpin{3}$. Because
$\mSpin{3}$ acts directly on $\reals^3$, the representation
of rotations in $\mSO{3}$ by elements of $\mSpin{3}$
may be viewed as more natural than the 
representation by unit quaternions.
The group $\mSU{2}\times \mSU{2}$ turns out to be isomorphic to
the spinor group  $\mSpin{4}$, but this isomorphism is less
obvious.

\medskip
In summary, we are going to define a group
$\mSpin{n}$ representing the rotations in $\mSO{n}$,
for any $n \geq 1$,
in the sense that there is a linear action of $\mSpin{n}$
on $\reals^n$ which
induces a surjective homomorphism
$\mapdef{\rho}{\mSpin{n}}{\mSO{n}}$ whose kernel
is $\{-1, 1\}$. 
Furthermore, the action of $\mSpin{n}$ on
$\reals^n$ is given in terms of multiplication in
an algebra $\mCl_n$ containing $\mSpin{n}$, and in 
which $\reals^n$ is also embedded.

\medskip
It turns out that as a bonus, for $n \geq 3$, the
group $\mSpin{n}$ is topologically simpler
than $\mSO{n}$, since  $\mSpin{n}$ is simply connected,
but $\mSO{n}$ is not.
By being astute, we can also construct a 
group $\mPin{n}$ and a linear action of 
$\mPin{n}$ on $\reals^n$ that induces
a surjective homomorphism
$\mapdef{\rho}{\mPin{n}}{\mO{n}}$ whose kernel
is $\{-1, 1\}$. The difficulty here is the presence
of the negative sign in (2). We will see how
Atiyah, Bott and Shapiro circumvent this problem
by using a ``twisted adjoint action,'' as opposed
to the usual adjoint action 
(where $v \mapsto u v u^{-1}$).

\medskip
Our presentation is heavily influenced by
Br\"ocker and tom Dieck \cite{Brocker} (Chapter 1, Section 6),
where most details can be found. This Chapter
is almost entirely taken from the first
11 pages of the beautiful
and seminal paper by Atiyah, Bott and Shapiro
\cite{AtiyahBottShap}, Clifford Modules, 
and we highly recommend it.
Another excellent (but concise) exposition can be 
found in Kirillov \cite{Kirillov01}.
A very thorough exposition can be found in two places:
\begin{enumerate}
\item
Lawson and Michelsohn \cite{Lawson},
where the material on $\mPin{p, q}$ and
$\mSpin{p, q}$ can be found in Chapter I.
\item
Lounesto's excellent book \cite{Lounesto}.
\end{enumerate}
One may also want to consult
Baker \cite{Baker},
Curtis \cite{Curtis}, Porteous \cite{Porteous},
Fulton and Harris (Lecture 20) \cite{Fulton91},
Choquet-Bruhat \cite{ChoquetBru2},
Bourbaki \cite{BourbakiA3},  
and Chevalley \cite{Chevalleyv2}, a classic.
The original source is Elie Cartan's book (1937)
whose translation in English appears in \cite{ECartan66}.

\medskip
We begin by recalling what is an algebra over a field.
Let $K$ denote any (commutative) field, although for our purposes
we may assume that $K = \reals$ (and occasionally,
$K = \complex$). Since we will only be dealing with
associative algebras with a multiplicative unit,
we only define algebras of this kind.

\begin{definition}
\label{algebradef2}
Given a field $K$, a {\it $K$-algebra\/}
is a $K$-vector space $A$ together with
a bilinear operation
$\mapdef{\cdot}{A\times A}{A}$, called {\it multiplication\/},
which makes $A$ into a ring with unity $1$ (or $1_A$, when we
want to be very precise).
This means that $\cdot$ is associative and that there is
a multiplicative identity element $1$ so that $1\cdot a = a\cdot 1 = a$,
for all $a\in A$.
Given two $K$-algebras $A$ and $B$, 
a {\it $K$-algebra homomorphism $\mapdef{h}{A}{B}$\/}
is a linear map that is also a ring homomorphism, with
$h(1_A) = 1_B$. 
\end{definition}

For example, the ring $M_n(K)$ of all $n\times n$ matrices over a 
field $K$ is a $K$-algebra.

\medskip
There is an obvious notion of {\it ideal\/} of a $K$-algebra: 
An ideal $\mfrac{A}\subseteq A$ is a linear subspace
of $A$ that is also a two-sided ideal with respect to multiplication
in $A$.
If the field $K$ is understood, we usually simply
say an algebra instead of a $K$-algebra.

\medskip
We will also need a quick review of  tensor products.
The basic idea is that tensor products
allow us to view multilinear maps as linear maps.
The maps become simpler, but the spaces 
(product spaces) become more complicated
(tensor products).
For more details, see 
Atiyah and Macdonald
\cite{AtiyahMac}.

\begin{definition}
\label{tensordef}
Given two $K$-vector spaces $E$ and $F$,
a {\it tensor product of $E$ and $F$\/}
is a pair $(E\tensor F,\, \tensor)$, where
$E\tensor F$ is a $K$-vector space and
$\mapdef{\tensor}{E\times F}{E\tensor F}$ is
a bilinear map, so that for every $K$-vector space
$G$ and every bilinear map
$\mapdef{f}{E\times F}{G}$, there is a unique linear
map $\mapdef{f_{\tensor}}{E\tensor F}{G}$ with
\[f(u, v) = f_{\tensor}(u\tensor v)\quad
\hbox{for all $u\in E$ and all $v\in V$},\]
as in the diagram below:
\[
\xymatrix{
E\times F \ar[r]^{\tensor} \ar[rd]_{f}  & \> E\tensor F \ar[d]^{f_{\tensor}} \\
       & \> G
}
\]
%
\end{definition}

\medskip
The vector space $E\tensor F$ is defined up to isomorphism.
The vectors $u\tensor v$, where $u\in E$ and $v\in F$,
generate $E\tensor F$. 

\remark
We should really denote the tensor product of $E$ and $F$ by 
$E\tensor_K F$, since it depends on the field $K$.
Since we usually deal with a fixed field $K$,
we use the simpler notation $E \tensor F$.

\medskip
We have  natural isomorphisms
\[(E \tensor F)\tensor G \approx E \tensor (F\tensor G)
\quad\hbox{and}\quad E\tensor F \approx F\tensor E.\]
Given two linear maps
$\mapdef{f}{E}{F}$ and $\mapdef{g}{E'}{F'}$, we have
a unique bilinear map \\
$\mapdef{f\times g}{E\times E'}{F\times F'}$ so that
\[(f\times g)(a, a') = (f(a), g(a'))
\quad\hbox{for all $a\in E$ and all $a'\in E'$}.\]
Thus, we have the bilinear map
$\mapdef{\tensor \circ (f\times g)}{E\times E'}{F\tensor F'}$, and
so, there is a unique linear map
$\mapdef{f\tensor g}{E\tensor E'}{F\tensor F'}$  so that
\[(f\tensor g)(a\tensor a') = f(a)\tensor g(a')
\quad\hbox{for all $a\in E$ and all $a'\in E'$}.\]

\medskip
Let us now assume that $E$ and $F$ are $K$-algebras.
We want to make $E\tensor F$ into a $K$-algebra.
Since the multiplication operations
$\mapdef{m_E}{E\times E}{E}$ and
$\mapdef{m_F}{F\times F}{F}$ are bilinear, we get
linear maps
$\mapdef{m_E'}{E\tensor E}{E}$ and
$\mapdef{m_F'}{F\tensor F}{F}$, and thus
the linear map
\[\mapdef{m_E'\tensor m_F'}{(E\tensor E)\tensor (F\tensor F)}
{E\tensor F}.\]
Using the isomorphism
$\mapdef{\tau}{(E\tensor E)\tensor (F\tensor F)}
{(E\tensor F)\tensor (E\tensor F)}$,
we get a linear map
\[\mapdef{m_{E\tensor F}}{(E\tensor F)\tensor (E\tensor F)}
{E\tensor F},\]
which defines a multiplication $m$ on $E\tensor F$
(namely, $m(u, v) = m_{E\tensor F}(u\tensor v)$).
It is easily checked that $E\tensor F$ is indeed
a $K$-algebra under the multiplication $m$.
Using the simpler notation $\cdot$ for $m$, we have
\[(a\tensor a')\cdot (b\tensor b') = (ab)\tensor (a'b')\]
for all $a, b\in E$ and all $a', b' \in F$.

\medskip
Given any vector space $V$ over a field $K$, there is 
a special $K$-algebra $T(V)$ together with
a linear map $\mapdef{i}{V}{T(V)}$,
with the following
universal mapping property: Given any $K$-algebra $A$,
for any linear map $\mapdef{f}{V}{A}$, there is
a unique $K$-algebra homomorphism $\mapdef{\overline{f}}{T(V)}{A}$
so that
\[f = \overline{f}\circ i,\]
as in the diagram below:
\[
\xymatrix{
V \ar[r]^{i} \ar[rd]_{f}  & \> T(V) \ar[d]^{\overline{f}} \\
       & \> A
}
\]
%
The algebra $T(V)$ is the {\it tensor algebra of $V$\/}.
The algebra $T(V)$ may be constructed as the direct sum
\[T(V) = \bigoplus_{i \geq 0}\, V^{\tensor i},\]
where $V^0 = K$, and $V^{\tensor i}$ is the $i$-fold
tensor product of $V$ with itself.
For every $i\geq 0$, there is a natural injection
$\mapdef{\iota_n}{V^{\tensor n}}{T(V)}$, and 
in particular, an injection
$\mapdef{\iota_0}{K}{T(V)}$. 
The multiplicative unit $\qone$ of $T(V)$ is the image
$\iota_0(1)$ in $T(V)$ of the unit $1$ of
the field $K$. 
Since every $v\in T(V)$
can be expressed as a finite sum
\[v = v_1 + \cdots + v_k,\]
where $v_i \in V^{\tensor n_i}$ and the
$n_i$ are natural numbers with $n_i \not= n_j$
if $i\not= j$, to define multiplication in $T(V)$,
using bilinearity, it is enough to define the multiplication
$V^{\tensor m}\times V^{\tensor n} \longrightarrow
V^{\tensor (m+n)}$. Of course, this is defined by
\[
(v_1\tensor \cdots \tensor v_m)\cdot
(w_1\tensor \cdots \tensor w_n) =
v_1\tensor \cdots \tensor v_m \tensor
w_1\tensor \cdots \tensor w_n.
\]
(This has to be made rigorous by using isomorphisms
involving the associativity
of tensor products; for details, see 
 see Atiyah and Macdonald
\cite{AtiyahMac}.)
The algebra $T(V)$ is an example of a {\it graded algebra\/},
where the {\it homogeneous elements of rank $n$\/} are the
elements in $V^{\tensor n}$.

\remark
It is important to note that multiplication in $T(V)$
is {\bf not\/} commutative. Also, in all rigor, the
unit $\qone$ of $T(V)$ is {\bf not equal\/} to $1$,
the unit of the field $K$. However, in view of the injection
$\mapdef{\iota_0}{K}{T(V)}$, for the sake of notational
simplicity, we will denote $\qone$ by $1$.
More generally,  in view of the injections
$\mapdef{\iota_n}{V^{\tensor n}}{T(V)}$, we identify elements of
$V^{\tensor n}$ with their images in $T(V)$.

\medskip
Most algebras of interest arise as well-chosen quotients
of the tensor algebra $T(V)$. This is true for the
{\it exterior algebra $\bigwedge^{\minibull} V$\/} (also called
{\it Grassmann algebra\/}), where we take the quotient
of $T(V)$ modulo the ideal generated by all elements
of the form $v\tensor v$, where $v\in V$, and for the 
{\it symmetric algebra $\mathrm{Sym}\> V$\/},
where we take the quotient
of $T(V)$ modulo the ideal generated by all elements
of the form $v\tensor w - w\tensor v$, where $v, w\in V$.

\medskip
A Clifford algebra may be viewed as a refinement
of the exterior algebra, in which 
we take the quotient
of $T(V)$ modulo the ideal generated by all elements
of the form $v\tensor v - \Phi(v)\cdot 1$, 
where $\Phi$ is the quadratic form associated
with a symmetric bilinear form 
$\mapdef{\varphi}{V\times V}{K}$, and $\mapdef{\cdot}{K\times T(V)}{T(V)}$
denotes the scalar product of the algebra $T(V)$.
For simplicity, let us assume that we are now
dealing with real algebras.

\begin{definition}
\label{cliffdef}
Let $V$ be a real finite-dimensional vector space
together with a symmetric bilinear form
$\mapdef{\varphi}{V\times V}{\reals}$ and associated
quadratic form $\Phi(v) = \varphi(v, v)$.
A {\it Clifford algebra associated with
$V$ and $\Phi$\/} is a real algebra $\mCl(V, \Phi)$
together with a linear map $\mapdef{i_{\Phi}}{V}{\mCl(V, \Phi)}$
satisfying the condition $(i_{\Phi}(v))^2 = \Phi(v)\cdot 1$
for all $v\in V$, and 
so that for every real algebra $A$ and every linear
map $\mapdef{f}{V}{A}$ with 
\[(f(v))^2 = \Phi(v)\cdot 1\quad\hbox{for all $v\in V$},\]
there is a unique algebra homomorphism
$\mapdef{\overline{f}}{\mCl(V, \Phi)}{A}$ so that
\[f = \overline{f}\circ i_{\Phi},\]
as in the diagram below:
\[
\xymatrix{
V \ar[r]^{i_{\Phi}} \ar[rd]_{f}  & \> \mCl(V, \Phi) \ar[d]^{\overline{f}} \\
       & \> A
}
\]
%
We use the notation $\lambda\cdot u$  for the product 
of a scalar $\lambda\in \reals$ 
and of an element $u$ in the algebra  $\mCl(V, \Phi)$,
and juxtaposition $uv$ for the multiplication of
two elements $u$ and $v$ in the algebra $\mCl(V, \Phi)$.
\end{definition}

\medskip
By a familiar argument, any two Clifford algebras
associated with $V$ and $\Phi$ are isomorphic.
We often denote $i_{\Phi}$ by $i$.

\medskip
To show the existence of $\mCl(V, \Phi)$, observe that
$T(V)/\mfrac{A}$ does the job, where $\mfrac{A}$ is the ideal of $T(V)$ 
generated by all elements of the form $v\tensor v - \Phi(v)\cdot 1$, 
where $v\in V$. The map $\mapdef{i_{\Phi}}{V}{\mCl(V, \Phi)}$ is the
composition
\[V \stackrel{\iota_1}{\longrightarrow} T(V)
\stackrel{\pi}{\longrightarrow} T(V)/\mfrac{A},\]
where $\pi$ is the natural quotient map.
We often denote the Clifford algebra
$\mCl(V, \Phi)$ simply by $\mCl(\Phi)$.

\remark
Observe that Definition \ref{cliffdef} does not
assert that $i_{\Phi}$ is injective or that there is
an injection of $\reals$ into $\mCl(V, \Phi)$, 
but we will prove later that both facts are true
when $V$ is finite-dimensional.
Also, as in the case of the tensor algebra, 
the unit of the algebra $\mCl(V, \Phi)$ and the unit of 
the field $\reals$ are {\bf not equal\/}.

\medskip
Since
\[\Phi(u + v) - \Phi(u) - \Phi(v) = 2\varphi(u, v)\]
and
\[(i(u + v))^2  =  (i(u))^2 + (i(v))^2 + i(u)  i(v) + i(v)  i(u),\]
using the fact that
\[i(u)^2 = \Phi(u)\cdot 1,\]
we get
\[i(u)  i(v) + i(v)  i(u) = 2\varphi(u, v)\cdot 1.\]
As a consequence, if $(u_1, \ldots, u_n)$
is an orthogonal basis w.r.t. $\varphi$
(which means that \\
$\varphi(u_j, u_k) = 0$
for all $j\not= k$), we have
\[i(u_j)  i(u_k) + i(u_k)  i(u_j) = 0\quad
\hbox{for all $j\not= k$.}\]

\remark
Certain authors drop the
unit $1$ of the Clifford algebra $\mCl(V, \Phi)$
when writing the identities
\[i(u)^2 = \Phi(u)\cdot 1\]
and
\[2\varphi(u, v)\cdot 1 =
i(u)  i(v) + i(v)  i(u),\]
where the second identity is often written as
\[\varphi(u, v) =
\frac{1}{2}(i(u)  i(v) + i(v)  i(u)).\]
This is very confusing and
technically wrong, because we only have an injection
of $\reals$ into $\mCl(V, \Phi)$, but $\reals$ is
{\bf not\/} a subset of $\mCl(V, \Phi)$.

\danger
We warn the readers that 
Lawson and Michelsohn \cite{Lawson} adopt 
the opposite of our sign convention in defining 
Clifford algebras, {\it i.e.}, they use the condition
\[(f(v))^2 = -\Phi(v)\cdot 1\quad\hbox{for all $v\in V$}.\]
The most confusing consequence of this is that
their $\mCl(p, q)$ is our $\mCl(q, p)$.

\medskip
Observe that when $\Phi \equiv 0$ is the quadratic
form identically zero everywhere, then the Clifford
algebra $\mCl(V, 0)$ is just the exterior algebra
$\bigwedge^{\minibull} V$.

\begin{example}
\label{ex1}
Let $V = \reals$,  $e_1 = 1$,
and assume that $\Phi(x_1e_1) = -x_1^2$.
Then, $\mCl(\Phi)$ is spanned by the basis $(1, e_1)$.
We have
\[e_1^2 = -1.\]
Under the bijection
\[e_1 \mapsto i,\]
the Clifford algebra $\mCl(\Phi)$, also
denoted by $\mCl_1$, is
isomorphic to the algebra of complex numbers $\complex$.

\medskip
Now, let $V = \reals^2$,
$(e_1, e_2)$ be the canonical basis, and
assume that $\Phi(x_1e_1 + x_2e_2) = -(x_1^2 + x_2^2)$.
Then, $\mCl(\Phi)$ is spanned by the basis
$(1, e_1, e_2, e_1e_2)$.
Furthermore, we have
\[e_2e_1 = - e_1e_2,\quad 
e_1^2 = -1,\quad e_2^2 = -1,\quad
(e_1e_2)^2 = -1.\]
Under the bijection
\[
e_1 \mapsto \qi,\quad
e_2 \mapsto \qj,\quad
e_1e_2 \mapsto \qk,
\]
it is easily checked that the quaternion
identities
\begin{align*}
 \qi^2 &= \qj^2 = \qk^2 =  -\qone,\\
 \qi\qj &= -\qj\qi = \qk,\\
 \qj\qk &= -\qk\qj = \qi,\\
 \qk\qi &= -\qi\qk = \qj,
\end{align*}
hold, and thus the Clifford algebra $\mCl(\Phi)$, also
denoted by $\mCl_2$, is
isomorphic to the algebra of quaternions $\quat$.
\end{example}

Our prime goal is to define an action of $\mCl(\Phi)$
on $V$ in such a way that by restricting this action
to some suitably chosen multiplicative subgroups
of $\mCl(\Phi)$, we get surjective homomorphisms
onto $\mO{\Phi}$ and $\mSO{\Phi}$, respectively.
The key point is that a reflection in $V$
about a hyperplane $H$ orthogonal to
a vector $w$ can be defined by such an action,
but some negative sign shows up.
A correct handling of signs is a bit subtle and 
requires the introduction of a canonical
anti-automorphism $t$,  and of a canonical automorphism
$\alpha$, defined as follows:

\begin{proposition}
\label{prop1}
Every Clifford algebra $\mCl(\Phi)$ possesses a canonical 
{\it anti-automorphism $\mapdef{t}{\mCl(\Phi)}{\mCl(\Phi)}$\/}
satisfying the properties
\[t(x  y) = t(y)  t(x),\quad
t\circ t = \id,\quad \hbox{and} \quad 
t(i(v)) = i(v),\] 
for all $x, y\in \mCl(\Phi)$ and all $v\in V$.
Furthermore, such an anti-automorphism is unique.
\end{proposition}

\begin{proof}
Consider the opposite algebra $\mCl(\Phi)^o$,
in which the product of $x$ and $y$ is given by $y x$.
It has the universal mapping property. Thus,
we get a unique isomorphism $t$, as in the diagram below:
\[
\xymatrix{
V \ar[r]^{i} \ar[rd]_{i}  & \> \mCl(V, \Phi) \ar[d]^{t} \\
       & \> \mCl(\Phi)^o
}
\]
%
\end{proof}

\medskip
We also denote $t(x)$ by $x^t$.
When $V$ is finite-dimensional,
for a more palatable description of $t$ in terms
of a basis of $V$,
see the paragraph following Theorem \ref{thm1}.

\medskip
The canonical automorphism
$\alpha$ is defined  using the proposition

\begin{proposition}
\label{prop2}
Every Clifford algebra $\mCl(\Phi)$ has a unique canonical 
{\it automorphism \\
$\mapdef{\alpha}{\mCl(\Phi)}{\mCl(\Phi)}$\/}
satisfying the properties
\[\alpha\circ \alpha = \id,\quad \hbox{and} \quad 
\alpha(i(v)) = -i(v),\] 
for  all $v\in V$.
\end{proposition}

\begin{proof}
Consider the linear map
$\mapdef{\alpha_0}{V}{\mCl(\Phi)}$
defined by $\alpha_0(v) = -i(v)$,  for
all $v\in V$. 
We get a unique homomorphism $\alpha$ as in the diagram below:
\[
\xymatrix{
V \ar[r]^{i} \ar[rd]_{\alpha_0}  & \> \mCl(V, \Phi) \ar[d]^{\alpha} \\
       & \> \mCl(\Phi)
}
\]
%
Furthermore, every $x\in \mCl(\Phi)$ can be written as
\[x = x_1\cdots x_m,\]
with $x_j\in i(V)$, and since $\alpha(x_j) = - x_j$,
we get $\alpha\circ \alpha = \id$. It is clear
that $\alpha$ is bijective.
\end{proof}

\medskip
Again, when $V$ is finite-dimensional,
a more palatable description of $\alpha$ in terms
of a basis of $V$ can be given.
If $(e_1, \ldots, e_n)$ is a basis of $V$,
then the Clifford algebra $\mCl(\Phi)$ consists
of certain kinds of ``polynomials,'' linear combinations
of monomials of the form 
$\sum_J \lambda_J e_J$, where 
$J = \{i_1, i_2, \ldots, i_k\}$
is any subset (possibly empty) of $\{1, \ldots, n\}$,
with $1 \leq i_1 < i_2 \cdots < i_k \leq n$,
and the monomial $e_J$ is the ``product''
$e_{i_1} e_{i_2} \cdots e_{i_k}$.
The map $\alpha$ is the linear map defined on
monomials by
\[\alpha(e_{i_1} e_{i_2} \cdots e_{i_k})  = 
(-1)^k e_{i_1} e_{i_2} \cdots e_{i_k}.\]
For a more rigorous explanation,
see the paragraph following Theorem \ref{thm1}.

\medskip
We now show that if $V$ has dimension $n$, then
$i$ is injective and
$\mCl(\Phi)$ has dimension $2^n$.
A clever way of doing this is to introduce a graded
tensor product.

\medskip
First, observe that
\[\mCl(\Phi) = \mCl^0(\Phi) \oplus \mCl^1(\Phi),\]
where 
\[\mCl^i(\Phi) = \{x\in \mCl(\Phi)\mid \alpha(x) = (-1)^i x\},\quad
\hbox{where $i = 0, 1$}.\]
We say that we have a {\it $\integs/2$-grading\/},
which means that if
$x\in \mCl^i(\Phi)$ and $y\in \mCl^j(\Phi)$, then
$x  y \in \mCl^{i + j\>(\mathrm{mod}\> 2)}(\Phi)$.

\medskip
When $V$ is finite-dimensional, since
every element of $\mCl(\Phi)$ is a linear
combination of the form $\sum_J \lambda_J e_J$ 
as explained earlier, in view of the description
of $\alpha$ given above, we see that
the elements of $\mCl^0(\Phi)$ are those
for which the monomials $e_J$ are products of
an even number of factors, and 
the elements of $\mCl^1(\Phi)$ are those
for which the monomials $e_J$ are products of
an odd number of factors.

\remark
Observe that $\mCl^0(\Phi)$ is a subalgebra of $\mCl(\Phi)$,
whereas $\mCl^1(\Phi)$ is not.

\medskip
Given two $\integs/2$-graded algebras
$A = A^0\oplus A^1$ and $B = B^0\oplus B^1$,
their {\it graded tensor product\/}
$A\>\widehat{\tensor}\> B$ is defined by
\begin{eqnarray*}
(A\>\widehat{\tensor}\> B)^0 & = &
(A^0\tensor B^0)\oplus (A^1\tensor B^1),\\
(A\>\widehat{\tensor}\> B)^1 & = &
(A^0\tensor B^1)\oplus (A^1\tensor B^0),
\end{eqnarray*}
with multiplication
\[(a'\tensor b)  (a\tensor b') = 
(-1)^{ij}(a'  a)\tensor (b  b'),\]
for $a\in A^i$ and $b\in B^j$.
The reader should check that 
$A\>\widehat{\tensor}\> B$ is indeed $\integs/2$-graded.

\begin{proposition}
\label{prop3}
Let $V$ and $W$ be finite dimensional
vector spaces with quadratic forms
$\Phi$ and $\Psi$. Then, there is a quadratic form
$\Phi\oplus \Psi$ on $V\oplus W$ defined by
\[(\Phi + \Psi)(v, w) = \Phi(v) + \Psi(w).\]
If we write $\mapdef{i}{V}{\mCl(\Phi)}$ and
$\mapdef{j}{W}{\mCl(\Psi)}$, we can define a linear map
\[\mapdef{f}{V\oplus W}{\mCl(\Phi)\>\widehat{\tensor}\> \mCl(\Psi)}\]
by
\[f(v, w) = i(v)\tensor 1 + 1 \tensor j(w).\]
Furthermore, the map $f$ induces an isomorphism
(also denoted by $f$)
\[\mapdef{f}{\mCl(V\oplus W)}{\mCl(\Phi)\>\widehat{\tensor}\> \mCl(\Psi).}\]
\end{proposition}

\begin{proof}
See Br\"ocker and tom Dieck \cite{Brocker}, Chapter 1, Section 6,
page 57.
\end{proof}

\medskip
As a corollary, we obtain the following
result:

\begin{theorem}
\label{thm1}
For every vector space $V$ of finite dimension $n$,
the map $\mapdef{i}{V}{\mCl(\Phi)}$ is injective.
Given a basis $(e_1, \ldots, e_n)$ of $V$, the $2^n - 1$
products
\[
i(e_{i_1}) i(e_{i_2}) \cdots i(e_{i_k}),\quad
1 \leq i_1 < i_2 \cdots < i_k \leq n,
\]
and $1$ form a basis of $\mCl(\Phi)$.
Thus, $\mCl(\Phi)$ has dimension $2^n$.
\end{theorem}

\begin{proof}
The proof is by induction on $n = \dimm(V)$.
For $n = 1$, the tensor algebra $T(V)$ is just 
the polynomial ring $\reals[X]$, where $i(e_1) = X$.
Thus, $\mCl(\Phi) = \reals[X]/(X^2 - \Phi(e_1))$, and
the result is obvious.
Since
\[i(e_j)  i(e_k) + i(e_k)  i(e_j) = 2\varphi(e_i, e_j)\cdot 1,\]
it is clear that the products
\[
i(e_{i_1}) i(e_{i_2}) \cdots i(e_{i_k}),\quad
1 \leq i_1 < i_2 \cdots < i_k \leq n,
\]
and $1$ generate $\mCl(\Phi)$. 
Now, there is always a basis that is orthogonal
with respect to $\varphi$ (for example, 
see Artin \cite{Artin91}, Chapter 7, or 
Gallier  \cite{Gallbook2}, Chapter 6, Problem 6.14),
and thus,
we have a splitting
\[(V, \Phi) = \bigoplus_{k = 1}^n (V_k, \Phi_k),\]
where $V_k$ has dimension $1$. Choosing a basis
so that $e_k\in V_k$, the theorem follows
by induction from Proposition \ref{prop3}.  
\end{proof}

\medskip
Since $i$ is injective,
for simplicity of notation, from now on we write $u$ for
$i(u)$.
Theorem \ref{thm1} implies that if 
$(e_1, \ldots, e_n)$ is an orthogonal basis of $V$, then
$\mCl(\Phi)$ is the algebra presented by the generators
$(e_1, \ldots, e_n)$ and the relations
\begin{eqnarray*}
e_j^2 & = & \Phi(e_j)\cdot 1,\quad 1\leq j \leq n,\quad\hbox{and}\\
e_je_k & = & -e_ke_j,\quad 1\leq j, k\leq n,\>j \not= k.
\end{eqnarray*}

\medskip
If $V$ has finite dimension $n$ and 
$(e_1, \ldots, e_n)$ is a basis of $V$,  
by Theorem \ref{thm1}, the maps $t$ 
and $\alpha$ are completely determined
by their action on the basis elements. 
Namely, $t$ is defined by
\begin{eqnarray*}
t(e_i) & = & e_i \\
t(e_{i_1} e_{i_2} \cdots e_{i_k}) & = &
e_{i_k} e_{i_{k-1}} \cdots e_{i_1},
\end{eqnarray*}
where $1 \leq i_1 < i_2 \cdots < i_k \leq n$,
and of course, $t(1) = 1$.
The map $\alpha$ is defined by
\begin{eqnarray*}
\alpha(e_i) & = & -e_i \\
\alpha(e_{i_1} e_{i_2} \cdots e_{i_k}) & = &
(-1)^k e_{i_1} e_{i_2} \cdots e_{i_k}
\end{eqnarray*}
where $1 \leq i_1 < i_2 \cdots < i_k \leq n$,
and of course, $\alpha(1) = 1$.
Furthermore, the even-graded elements 
(the elements of  $\mCl^0(\Phi)$)
are those generated by $1$ and the basis
elements consisting of an even number of factors
$e_{i_1} e_{i_2} \cdots e_{i_{2k}}$, and
the odd-graded elements (the elements of $\mCl^1(\Phi)$)
are those generated by the basis
elements consisting of an odd number of factors
$e_{i_1} e_{i_2} \cdots e_{i_{2k + 1}}$.

\medskip
We are now ready to define the Clifford group and
investigate some of its properties.

\section{Clifford Groups}
\label{seccliff2}
First, we define {\it conjugation\/} on a Clifford
algebra $\mCl(\Phi)$ as the map
\[x \mapsto \overline{x} = t(\alpha(x))\quad
\hbox{for all $x\in \mCl(\Phi)$}.\]
Observe that
\[t\circ \alpha = \alpha \circ t.\]
If $V$ has finite dimension $n$ and 
$(e_1, \ldots, e_n)$ is a basis of $V$,  
in view of previous remarks, conjugation is defined by
\begin{eqnarray*}
\overline{e_i} & = & -e_i \\
\overline{e_{i_1} e_{i_2} \cdots e_{i_k}} & = &
(-1)^k e_{i_k} e_{i_{k-1}} \cdots e_{i_1}
\end{eqnarray*}
where $1 \leq i_1 < i_2 \cdots < i_k \leq n$,
and of course, $\overline{1} = 1$.
Conjugation is an anti-automorphism.

\medskip
The multiplicative group of invertible elements of $\mCl(\Phi)$
is denoted by $\mCl(\Phi)^*$.
Observe that for any $x\in V$, if 
$\Phi(x)\not= 0$,  then
$x$ is invertible because  $x^2 = \Phi(x)$; that is, 
$x\in \mCl(\Phi)^*$.

\medskip
We would like   $\mCl(\Phi)^*$ to act on $V$ {\it via\/}
\[x\cdot v = \alpha(x) v x^{-1},\]
where $x\in  \mCl(\Phi)^*$ and $v\in V$.
In general, there is no reason why $\alpha(x) v x^{-1}$ should
be in $V$ or why this action defines an automorphism of $V$,
so we restrict this map to the subset $\Gamma(\Phi)$
of $\mCl(\Phi)^*$ 
as follows.

\begin{definition}
\label{cliffgr}
Given a finite dimensional vector space $V$ and a quadratic form $\Phi$
on $V$, the {\it Clifford group of $\Phi$\/}
is the group
\[\Gamma(\Phi) = \{x\in \mCl(\Phi)^*\mid
\alpha(x)  v   x^{-1} \in V\quad
\hbox{for all $v\in V$}\}.\]
The map $\mapdef{N}{\mCl(Q)}{\mCl(Q)}$ given by
\[N(x) = x  \overline{x}\]
is called the {\it norm\/} of $\mCl(\Phi)$.
\end{definition}

For any $x\in \Gamma(\Phi)$,
let $\mapdef{\rho_x}{V}{V}$ be the map 
defined by
\[v \mapsto \alpha(x)  v   x^{-1}, \quad v \in V.\]
It is not entirely obvious why 
the map  $\mapdef{\rho}{\Gamma(\Phi)}{\mGL{V}}$ given by
$x \mapsto \rho_x$
is a linear action, and for that matter, why
$\Gamma(\Phi)$ is a group. 
This is because
$V$ is finite-dimensional and $\alpha$ is an automorphism.

\begin{proof}
For any $x\in \Gamma(\Phi)$,
the map $\rho_x$ from $V$ to $V$ defined by
\[v \mapsto \alpha(x)  v   x^{-1}\]
is clearly linear. If 
$ \alpha(x)  v   x^{-1} = 0$, since by hypothesis
$x$ is invertible and since
$\alpha$ is an automorphism $\alpha(x)$ is also invertible,
so $v = 0$. Thus our linear map is injective, and
since $V$ has finite dimension, it is bijective.
To prove that  $x^{-1}\in \Gamma(\Phi)$, pick any $v\in V$.
Since the linear map $\rho_x$ is bijective, there is some $w\in V$
such that $\rho_x(w) = v$, which means that
$\alpha(x)  w   x^{-1} = v$. Since $x$ is invertible and $\alpha$ is
an automorphism, we get
\[
\alpha(x^{-1}) v x= w,
\]
so $\alpha(x^{-1}) v x\in V$; since this holds for any $v\in V$, we have
$x^{-1} \in \Gamma(\Phi)$.
Since $\alpha$ is an automorphism, if $x, y\in  \Gamma(\Phi)$,
for any $v\in V$ we have
\[
\rho_y(\rho_x(v)) = \alpha(y)\alpha(x) v x^{-1} y^{-1} =
\alpha(yx) v (yx)^{-1} = \rho_{yx}(v),
\] 
which shows that
$\rho_{yx}$ is a linear automorphism of $V$,  so
$yx\in \Gamma(\Phi)$ and  $\rho$ is a homomorphism.
Therefore, $\Gamma(\Phi)$ is
a group and $\rho$ is a linear representation.
\end{proof}

\medskip
We also define the group $\Gamma^+(\Phi)$, called
the {\it special Clifford group\/}, by
\[\Gamma^+(\Phi) = \Gamma(\Phi) \cap \mCl^0(\Phi).\]
Observe that $N(v) = - \Phi(v)\cdot 1$
for all $v\in V$. 
Also, if $(e_1, \ldots, e_n)$ is a basis of $V$,
we leave it as an exercise to check that
\[N(e_{i_1} e_{i_2} \cdots e_{i_k}) =
(-1)^k \Phi(e_{i_1}) \Phi(e_{i_2}) \cdots \Phi(e_{i_k})\cdot 1.\]

\remark
The map $\mapdef{\rho}{\Gamma(\Phi)}{\mGL{V}}$ given by
$x \mapsto \rho_x$
is called the {\it twisted adjoint representation\/}.
It was introduced by Atiyah, Bott and Shapiro
\cite{AtiyahBottShap}. It has the advantage of not 
introducing a spurious negative sign, {\it i.e.},
when $v\in V$ and $\Phi(v) \not= 0$, the map $\rho_v$ is
the reflection $s_v$ about the hyperplane orthogonal to $v$
(see Proposition \ref{refleclem}). Furthermore,
when $\Phi$ is nondegenerate, the kernel $\Ker(\rho)$
of the representation $\rho$ is given by
$\Ker(\rho) = \reals^*\cdot 1$,
where $\reals^* = \reals - \{0\}$. 
The earlier
{\it adjoint representation\/} (used by Chevalley
\cite{Chevalleyv2} and others) is given by
\[v \mapsto x  v   x^{-1}.\]
Unfortunately, in this case, $\rho_x$ represents $-s_v$,
where $s_v$ is the reflection about the hyperplane orthogonal to $v$.
Furthermore, the kernel of the representation $\rho$ 
is generally bigger than $\reals^*\cdot 1$.
This is the reason why the
twisted adjoint representation is preferred
(and  must be used for a proper treatment of the 
$\Pin$  group).

\begin{proposition}
\label{prop4}
The maps $\alpha$ and $t$ induce an automorphism and
an anti-automorphism of the Clifford group, $\Gamma(\Phi)$.
\end{proposition}

\begin{proof}
It is not very instructive; see 
Br\"ocker and tom Dieck \cite{Brocker}, Chapter 1, Section 6, page 58.
\end{proof}

\medskip
The following proposition shows why
Clifford groups generalize the quaternions.

\begin{proposition}
\label{refleclem}
Let $V$ be a finite dimensional vector
space and $\Phi$ a quadratic form on $V$.
For every element $x$ of the Clifford group
$\Gamma(\Phi)$, if $x\in V$ and $\Phi(x)\not= 0$, then
the map $\mapdef{\rho_x}{V}{V}$ given by
\[v \mapsto \alpha(x)vx^{-1}\quad
\hbox{for all $v\in V$}
\] 
is the reflection
about the hyperplane $H$ orthogonal to
the vector $x$.
\end{proposition}

\begin{proof}
Recall that the reflection $s$ 
about the hyperplane $H$ orthogonal to
the vector $x$ is given by
\[s(\novect{u}) = \novect{u} - 
2\,\frac{\varphi(u, x)}{\Phi(x)}\cdot x.\]
However, we have
\[x^2 = \Phi(x)\cdot 1
\quad\hbox{and}\quad
u  x + x  u= 2\varphi(u, x)\cdot 1.\]
Thus, we have
 \begin{eqnarray*} 
s(u) &= &  u - 2\,\frac{\varphi(u, x)}{\Phi(x)}\cdot x\\
& = & u - 2\varphi(u, x)\cdot\left(\frac{1}{\Phi(x)}\cdot x\right)\\
& = & u - 2\varphi(u, x) \cdot x^{-1}\\
& = & u - 2\varphi(u, x)\cdot (1 x^{-1})\\
& = & u - (2\varphi(u, x)\cdot 1) x^{-1}\\
& = & u  -(u  x + x  u) x^{-1}\\
& = & -x u x^{-1}\\
& = & \alpha(x) u x^{-1},
\end{eqnarray*} 
since $\alpha(x) = -x$, for $x\in V$.
\end{proof}

\medskip
Recall that the linear representation
\[\mapdef{\rho}{\Gamma(\Phi)}{\mGL{V}}\]
is given  by
\[\rho(x)(v) = \alpha(x)  v   x^{-1},\]
for all $x\in \Gamma(\Phi)$ and all $v\in V$.
We would like to show that $\rho$ is a surjective homomorphism
from $\Gamma(\Phi)$ onto $\mO{\varphi}$, and
a surjective homomorphism
from $\Gamma^+(\Phi)$ onto $\mSO{\varphi}$. 
For this, we will 
need to assume that $\varphi$ is nondegenerate,
which means that for every $v\in V$, if
$\varphi(v, w) = 0$ for all $w\in V$, then $v = 0$.
For simplicity of exposition, we first
assume that $\Phi$ is the quadratic form
on $\reals^n$ defined by
\[\Phi(x_1, \ldots, x_n) = -(x_1^2 + \cdots + x_n^2).\]
Let $\mCl_n$ denote the Clifford algebra $\mCl(\Phi)$ and
$\Gamma_n$ denote the Clifford group $\Gamma(\Phi)$.
The following lemma plays a crucial role:

\begin{lemma}
\label{kernelem}
The kernel of the map 
$\mapdef{\rho}{\Gamma_n}{\mGL{n}}$ is $\reals^*\cdot 1$,
the multiplicative group of nonzero scalar multiples of $1\in \mCl_n$.
\end{lemma}

\begin{proof}
If $\rho(x) = \id$, then
\begin{equation}
\alpha(x)  v = v  x\quad
\hbox{for all $v\in \reals^n$}.
\tag{1}
\end{equation}
Since $\mCl_n = \mCl_n^0\oplus \mCl_n^1$,
we can write $x = x^0 + x^1$, with $x^i\in \mCl_n^i$ for $i = 0, 1$.
Then, equation (1) becomes
\begin{equation}
x^0  v = v  x^0\quad
\hbox{and}\quad -x^1  v = v  x^1\quad \hbox{for all $v\in \reals^n$}.
\tag{2}
\end{equation}
Using Theorem \ref{thm1}, we can express $x^0$ as a linear combination
of monomials in the canonical basis $(e_1, \ldots, e_n)$, so that
\[x^0 = a^0 + e_1  b^1,
\quad\hbox{with $a^0\in \mCl_n^0,\> b^1\in \mCl_n^1$},\]
where neither $a^0$ nor $b^1$ contains a summand with a factor $e_1$.
Applying the first relation in (2) to $v = e_1$, we get
\begin{equation}
e_1 a^0 + e_1^2  b^1 = a^0 e_1 + e_1 b^1 e_1.
\tag{3}
\end{equation}
Now, the basis  $(e_1, \ldots, e_n)$ is orthogonal w.r.t. $\Phi$,
which implies that
\[e_je_k = -e_ke_j\quad
\hbox{for all $j\not= k$}.\]
Since each monomial in $a^0$ is of even degree and
contains no factor $e_1$, we get
\[a^0e_1 = e_1a^0.\]
Similarly, since $b^1$ is of odd degree and 
contains no factor $e_1$, we get
\[e_1b^1e_1 = -e_1^2b^1.\]
But then, from (3), we get
\[e_1a^0 + e_1^2b^1 = a^0e_1 + e_1b^1 e_1 = e_1a^0 - e_1^2b^1,\]
and so, $e_1^2b^1 = 0$. 
However, $e_1^2 = -1$, and so, $b_1 = 0$.
Therefore, $x_0$ contains
no monomial with a factor $e_1$.
We can apply the same argument to the other basis elements
$e_2, \ldots, e_n$, and thus, we just proved that
$x^0 \in \reals\cdot 1$.

\medskip
A similar argument applying to the second equation in (2),
with $x^1 = a^1 + e_1b^0$ and $v = e_1$
shows that $b^0 = 0$. We also conclude that $x^1\in \reals\cdot 1$.
However, $\reals\cdot 1 \subseteq \mCl_n^0$, and so $x^1 = 0$.
Finally, $x = x^0\in (\reals\cdot 1)\cap \Gamma_n = \reals^*\cdot 1$.
\end{proof}

\remark
If $\Phi$ is any nondegenerate
quadratic form, we know (for instance,
see Artin \cite{Artin91}, Chapter 7, or 
Gallier  \cite{Gallbook2}, Chapter 6, Problem 6.14)
that there is an orthogonal  basis 
$(e_1, \ldots, e_n)$ with respect to $\varphi$
({\it i.e.} $\varphi(e_j, e_k) = 0$ for all $j\not= k$).
Thus, the commutation relations
\begin{eqnarray*}
e_j^2 & = & \Phi(e_j)\cdot 1,\quad\hbox{with $\Phi(e_j) \not= 0$},\quad
1\leq j \leq n,\quad\hbox{and}\\
e_je_k & = & -e_ke_j,\quad 1\leq j, k\leq n,\>j \not= k 
\end{eqnarray*}
hold, and since the proof only rests on these facts,
Lemma \ref{kernelem} holds for any 
nondegenerate quadratic form.

\danger
However, Lemma \ref{kernelem} may fail for degenerate
quadratic forms. For example, if $\Phi \equiv 0$, then
$\mCl(V, 0) = \bigwedge^{\minibull} V$. Consider the element
$x = 1 + e_1e_2$. Clearly, $x^{-1} = 1 - e_1e_2$.
But now, for any $v\in V$, we have
\[\alpha(1 + e_1e_2)v(1 + e_1e_2)^{-1} = 
(1 + e_1e_2)v(1 - e_1e_2) = v.\]
Yet, $1 + e_1e_2$ is not a scalar multiple of $1$. 

\medskip
The following proposition shows that the notion of norm
is well-behaved.

\begin{proposition}
\label{prop5}
If $x\in \Gamma_n$, then $N(x) \in \reals^*\cdot 1$.
\end{proposition}

\begin{proof}
The trick is to show that $N(x)$ is in the kernel of $\rho$.
To say that $x\in \Gamma_n$ means that
\[\alpha(x)vx^{-1}\in \reals^n
\quad\hbox{for all $v\in \reals^n$}.\]
Applying $t$, we get
\[t(x)^{-1}v t(\alpha(x)) = \alpha(x)vx^{-1},\]
since $t$ is the identity on $\reals^n$.
Thus, we have
\[v = t(x)\alpha(x)v (t(\alpha(x)) x)^{-1} = 
\alpha(\overline{x} x)v (\overline{x} x)^{-1},\]
so $\overline{x}x\in \Ker(\rho)$.
By Proposition \ref{prop4}, we have $\overline{x}\in \Gamma_n$,
and so, $x\overline{x} = \overline{\overline{x}}\,\overline{x}
\in \Ker(\rho)$.
\end{proof}

\remark
Again, the proof also holds for the Clifford group 
$\Gamma(\Phi)$ associated with any nondegenerate quadratic form $\Phi$.
When $\Phi(v) = -\norme{v}^2$, where $\norme{v}$ is the
standard Euclidean norm of $v$, we have
$N(v) = \norme{v}^2\cdot 1$ for all $v\in V$.
However, for other quadratic forms, it is possible that
$N(x) = \lambda \cdot 1$ where $\lambda < 0$,
and this is a difficulty that needs to be overcome.

\begin{proposition}
\label{prop6}
The restriction of the norm $N$ to $\Gamma_n$
is a homomorphism $\mapdef{N}{\Gamma_n}{\reals^* \cdot 1}$,
and $N(\alpha(x)) = N(x)$ for all $x\in \Gamma_n$.
\end{proposition}

\begin{proof}
We have
\[N(xy) = xy\overline{y}\,\overline{x} = xN(y)\overline{x}
= x\overline{x} N(y) = N(x) N(y),\]
where the third equality holds because $N(x)\in \reals^*\cdot 1$.
We also have
\begin{equation*}
N(\alpha(x)) = \alpha(x)\alpha(\overline{x}) =
\alpha(x\overline{x}) = \alpha(N(x)) = N(x).
\qedhere
\end{equation*}
\end{proof}

\remark
The proof also holds for the Clifford group 
$\Gamma(\Phi)$ associated with any nondegenerate quadratic form $\Phi$.

\begin{proposition}
\label{prop7}
We have $\reals^n - \{0\} \subseteq \Gamma_n$ and
$\rho(\Gamma_n) \subseteq \mO{n}$.
\end{proposition}

\begin{proof}
Let $x\in \Gamma_n$ and $v\in \reals^n$, with $v\not= 0$.
We have
\[N(\rho(x)(v)) = N(\alpha(x)vx^{-1}) =
N(\alpha(x)) N(v) N(x^{-1}) = N(x) N(v) N(x)^{-1} = N(v),\]
since  $\mapdef{N}{\Gamma_n}{\reals^* \cdot 1}$.
However, for $v\in \reals^n$, we know that
\[N(v) = -\Phi(v)\cdot 1.\]
Thus, $\rho(x)$ is norm-preserving, and so,
$\rho(x)\in \mO{n}$.
\end{proof}

\remark
The proof that $\rho(\Gamma(\Phi)) \subseteq \mO{\Phi}$
also holds for the Clifford group 
$\Gamma(\Phi)$ associated with any nondegenerate quadratic form $\Phi$.
The first statement needs to be replaced by the fact
that every non-isotropic vector in $\reals^n$
(a vector is non-isotropic
if $\Phi(x) \not= 0$) belongs to $\Gamma(\Phi)$.
Indeed, $x^2 = \Phi(x)\cdot 1$, which implies that
$x$ is invertible.

\medskip
We are finally ready for the introduction of
the groups $\mPin{n}$ and $\mSpin{n}$.
\section{The Groups $\mPin{n}$ and $\mSpin{n}$}
\label{seccliff3}

\begin{definition}
\label{spindefin}
We define the {\it pinor group $\mPin{n}$\/} as the kernel
$\Ker(N)$ of the homomorphism
$\mapdef{N}{\Gamma_n}{\reals^*\cdot 1}$, and the 
{\it spinor group $\mSpin{n}$\/}
as $\mPin{n}\cap \Gamma^+_n$.
\end{definition}

\medskip
Observe that if $N(x) = 1$, then $x$ is invertible, and
$x^{-1} = \overline{x}$ since $x\overline{x} = N(x) = 1$. 
Thus, we can write
\begin{align*}
\mPin{n} & = \{x\in \mCl_{n} \mid
\alpha(x) v x^{-1} \in \reals^n\quad\hbox{for all $v\in \reals^n$,}\quad
N(x) = 1\} \\
& = \{x\in \mCl_{n} \mid
\alpha(x) v \overline{x} \in \reals^n\quad\hbox{for all $v\in \reals^n$,}\quad
x\overline{x} = 1\},
\end{align*}
and
\begin{align*}
\mSpin{n}& = \{x\in \mCl^0_{n} \mid
x v x^{-1} \in \reals^n\quad\hbox{for all $v\in \reals^n$,}\quad
N(x) = 1\} \\
& = \{x\in \mCl^0_{n} \mid
x v \overline{x} \in \reals^n\quad\hbox{for all $v\in \reals^n$,}\quad
x\overline{x} = 1\}
\end{align*}

\remark
According to Atiyah, Bott and Shapiro, 
the use of the name $\mPin{k}$ is a joke
due to Jean-Pierre Serre
(Atiyah, Bott and Shapiro \cite{AtiyahBottShap}, page 1).

\begin{theorem}
\label{spinthm}
The restriction of $\mapdef{\rho}{\Gamma_n}{\mO{n}}$ 
to the pinor group $\mPin{n}$ is a
surjective homomorphism
$\mapdef{\rho}{\mPin{n}}{\mO{n}}$
whose kernel is $\{-1, 1\}$, and 
the restriction of $\rho$ to the spinor group $\mSpin{n}$ is a
surjective homomorphism
$\mapdef{\rho}{\mSpin{n}}{\mSO{n}}$
whose kernel is $\{-1, 1\}$.
\end{theorem}

\begin{proof}
By Proposition \ref{prop7}, we have a map
$\mapdef{\rho}{\mPin{n}}{\mO{n}}$.
The reader can easily check that $\rho$ is a homomorphism.
By the Cartan-Dieudonn\'e theorem
(see Bourbaki \cite{BourbakiA3}, or 
Gallier \cite{Gallbook2}, Chapter 7, Theorem 7.2.1),
every isometry  $f\in \mSO{n}$ is the composition
$f = s_1\circ \cdots \circ s_k$ of hyperplane reflections $s_j$.
If we assume that $s_j$ is a reflection about the
hyperplane $H_j$ orthogonal to the nonzero vector $w_j$,
by Proposition \ref{refleclem}, 
$\rho(w_j) = s_j$. Since $N(w_j) = \norme{w_j}^2\cdot 1$,
we can replace $w_j$ by $w_j/\norme{w_j}$, so that
$N(w_1\cdots w_k) = 1$, and then
\[f = \rho(w_1\cdots w_k),\]
and $\rho$ is surjective.
Note that
\[\Ker(\rho\mid \mPin{n}) = \Ker(\rho) \cap \ker(N) =
\{t\in \reals^*\cdot 1 \mid N(t) = 1\} = \{-1, 1\}.\]
As to $\mSpin{n}$, we just need to show that
the restriction of $\rho$ to $\mSpin{n}$ 
maps $\Gamma_n$ into $\mSO{n}$.
If this was not the case,
there would be some improper isometry
$f\in \mO{n}$ so that
$\rho(x) = f$, where $x\in \Gamma_n \cap \mCl_n^0$.
However, we can express $f$ as the composition of
an odd number of reflections, say
\[f = \rho(w_1\cdots w_{2k+1}).\]
Since
\[\rho(w_1\cdots w_{2k+1}) = \rho(x),\]
we have
$x^{-1}w_1\cdots w_{2k+1} \in \Ker(\rho)$.
By Lemma \ref{kernelem}, we must have
\[x^{-1}w_1\cdots w_{2k+1} = \lambda\cdot 1\]
for some $\lambda\in \reals^*$, and thus
\[w_1\cdots w_{2k+1}  = \lambda\cdot x,\]
where $x$ has even degree and
$w_1\cdots w_{2k+1}$ has odd degree, 
which is impossible. 
\end{proof}

\medskip
Let us denote the set of elements $v\in \reals^n$ with $N(v) = 1$
(with norm $1$) by $S^{n-1}$.
We have the following corollary of 
Theorem \ref{spinthm}:

\begin{cor}
\label{spincor}
The group $\mPin{n}$ is generated by $S^{n-1}$, 
and every element of $\mSpin{n}$ can be written
as the product of an even number of elements
of $S^{n - 1}$.
\end{cor}

\begin{example}
\label{ex2}
The reader should verify that
\[\mPin{1} \approx \integs/4\integs,
\quad
\mSpin{1} = \{-1, 1\} \approx \integs/2\integs,\]
and also that
\[\mPin{2} \approx \{a e_1 + b e_2\mid a^2 + b^2 = 1\}
\cup \{c 1 + d e_1e_2\mid c^2 + d^2 = 1\}, \quad
\mSpin{2} = \mU{1}.\]
We may also write
$\mPin{2} = \mU{1} + \mU{1}$, 
where $\mU{1}$ is the group of complex numbers of
modulus $1$ (the unit circle in $\reals^2$).
It can also be shown that
$\mSpin{3} \approx \mSU{2}$ and
$\mSpin{4} \approx \mSU{2}\times \mSU{2}$.
The group $\mSpin{5}$ is isomorphic to the
symplectic group ${\bf Sp}(2)$, and
$\mSpin{6}$ is isomorphic to $\mSU{4}$
(see Curtis \cite{Curtis} or Porteous \cite{Porteous}).
\end{example}

\medskip
Let us take a closer look at $\mSpin{2}$.
The Clifford algebra $\mCl_2$ is generated by the four
elements
\[1,\> e_1,\> e_2,\>, e_1e_2,\]
and they satisfy the relations
\[e_1^2 = -1,\quad e_2^2 = -1, \quad  e_1e_2 = -e_2e_1.\]
The group  $\mSpin{2}$ consists of all products
\[\prod_{i = 1}^{2k} (a_i e_1 + b_i e_2)\]
consisting of an even number of factors
and such that $a_i^2 + b_i^2 = 1$.
In view of the above relations, every such element can be written as
\[x = a 1 + b e_1e_2,\]
where $x$ satisfies the conditions that $xvx^{-1}\in \reals^2$
for all $v\in \reals^2$, and  $N(x) = 1$.
Since 
\[\overline{X}  = a 1 - b e_1e_2,\]
we get
\[N(x) =  a^2 + b^2,\]
and the condition $N(x) = 1$ is simply
$a^2 + b^2  = 1$. 

\medskip
We claim that if $x\in \mCl^0_{2}$, then
$xvx^{-1}\in \reals^2$.
Indeed, since $x\in \mCl^0_{2}$ and $v\in \mCl^1_{2}$,
we have $xvx^{-1}\in \mCl^1_{2}$, which implies
that $xvx^{-1} \in \reals^2$, since the only elements
of  $\mCl^1_{2}$ are those in $\reals^2$.
Then, $\mSpin{2}$ consists of those elements
$x = a 1 + b  e_1e_2$ so that
$a^2 + b^2  = 1$.
If we let $\qi = e_1e_2$, we observe that
\begin{eqnarray*}
\qi^2 & = & -1, \\
e_1\qi & = & -\qi e_1 \> = \> -e_2, \\ 
e_2\qi & = & -\qi e_2 \> = \> e_1. 
\end{eqnarray*}
Thus, $\mSpin{2}$ is isomorphic to $\mU{1}$.
Also note that
\[e_1 (a 1 + b \qi) = (a 1 - b \qi) e_1.\] 
Let us find out explicitly what is the action
of $\mSpin{2}$ on $\reals^2$. Given
$X = a 1 + b \qi$, with $a^2 + b^2 = 1$,
for any $v = v_1 e_1 + v_2 e_2$,
we have
\begin{eqnarray*}
\alpha(X) v X^{-1} 
& = & X(v_1 e_1 + v_2 e_2) X^{-1}\\ 
& = & X(v_1 e_1 + v_2 e_2) (-e_1 e_1) \overline{X}\\
& = &  X(v_1 e_1 + v_2 e_2)(-e_1) (e_1 \overline{X})\\
& = &  X(v_1 1  + v_2 \qi) X e_1\\
& = & X^2(v_1 1  + v_2 \qi) e_1 \\
& = & (((a^2 - b^2) v_1 - 2ab v_2) 1 + (a^2 - b^2) v_2 + 2ab v_1)\qi) e_1 \\
& = & ((a^2 - b^2) v_1 - 2ab v_2) e_1 + (a^2 - b^2) v_2 + 2ab v_1)e_2. 
\end{eqnarray*}

Since $a^2 + b^2 = 1$, we can write 
$X = a 1 + b \qi =  (\cos\theta) 1 + (\sin\theta) \qi$, 
and the above derivation shows that
\[\alpha(X) v X^{-1} =
(\cos 2\theta v_1 - \sin 2\theta v_2)e_1 +
(\cos 2\theta v_2 + \sin 2\theta v_1)e_2.\] 
This means that the rotation $\rho_X$ induced by $X \in \mSpin{2}$
is the rotation of angle $2\theta$ around the origin.
Observe that the maps
\[v \mapsto v(- e_1),\quad X \mapsto X e_1\]
establish bijections between $\reals^2$ and
$\mSpin{2} \simeq \mU{1}$.
Also, note that the action of $X = \cos\theta + i \sin\theta$ 
viewed as a complex number
yields the rotation of angle $\theta$, whereas the action 
of $X = (\cos\theta) 1 +  (\sin\theta) \qi$ 
viewed as a member of $\mSpin{2}$ yields
the rotation of angle $2\theta$.  There is nothing wrong.
In general, $\mSpin{n}$ is a two--to--one cover of
$\mSO{n}$.

\medskip
Next, let us take a closer look at $\mSpin{3}$.
The Clifford algebra $\mCl_3$ is generated by the eight
elements
\[1,\> e_1,\> e_2,\>, e_3,\>, e_1e_2,\> e_2e_3,\> e_3e_1,\> e_1e_2e_3,\]
and they satisfy the relations
\[e_i^2 = -1,\quad e_je_j = -e_je_i,\quad 1\leq i, j\leq 3,\> i\not= j.\]
The group  $\mSpin{3}$ consists of all products
\[\prod_{i = 1}^{2k} (a_i e_1 + b_i e_2 + c_i e_3)\]
consisting of an even number of factors
and such that $a_i^2 + b_i^2 + c_i^2  = 1$.
In view of the above relations, every such element can be written as
\[x = a 1 + b e_2e_3 + c e_3e_1 + d e_1e_2,\]
where $x$ satisfies the conditions that $xvx^{-1}\in \reals^3$
for all $v\in \reals^3$, and  $N(x) = 1$.
Since 
\[\overline{X}  = a 1 - b e_2e_3 - c e_3e_1 - d e_1e_2,\]
we get
\[N(x) =  a^2 + b^2 + c^2 + d^2,\]
and the condition $N(x) = 1$ is simply
$a^2 + b^2 + c^2 + d^2 = 1$.

\medskip
It turns out that the conditions $x\in \mCl^0_{3}$
and $N(x) = 1$ imply
that  $xvx^{-1}\in \reals^3$ for all $v\in \reals^3$. 
To prove this, first observe
that $N(x) = 1$ implies that $x^{-1} = \pm\overline{x}$, and 
that $\overline{v} = -v$ for any $v\in \reals^3$,  and so,
\[\overline{xvx^{-1}} = -xvx^{-1}.\]
Also, since $x\in \mCl^0_{3}$ and $v\in \mCl^1_{3}$,
we have $xvx^{-1}\in \mCl^1_{3}$. Thus, we can write
\[xvx^{-1} = u + \lambda e_1e_2 e_3,
\quad\hbox{for some $u\in \reals^3$ and some $\lambda\in \reals$.} 
\]
But 
\[\overline{e_1e_2e_3} = -e_3e_2e_1 = e_1e_2e_3,\]
and so,
\[\overline{xvx^{-1}} = -u + \lambda e_1e_2e_3 = -xvx^{-1} = -u - \lambda e_1e_2e_3,\]
which implies that $\lambda = 0$.
Thus, $xvx^{-1}\in \reals^3$, as claimed.
Then, $\mSpin{3}$ consists of those elements
$x = a 1 + b e_2e_3 + c e_3e_1 + d e_1e_2$ so that
$a^2 + b^2 + c^2 + d^2 = 1$.
Under the bijection
\[\qi \mapsto e_2e_3,\> \qj \mapsto e_3e_1,\> \qk \mapsto e_1e_2,\]
we can check that we have an isomorphism between the group
$\mSU{2}$ of unit quaternions and $\mSpin{3}$.
If $X = a 1 + b e_2e_3 + c e_3e_1 + d e_1e_2\in \mSpin{3}$,
observe that
\[X^{-1} = \overline{X}  = a 1 - b e_2e_3 - c e_3e_1 - d e_1e_2.\]
Now, using the identification
\[\qi \mapsto e_2e_3,\> \qj \mapsto e_3e_1,\> \qk \mapsto e_1e_2,\]
we can easily check that
\begin{eqnarray*}
(e_1e_2e_3)^2 & = & 1,\\
(e_1e_2e_3) \qi &= & \qi (e_1e_2e_3)\> =\> -e_1,\\
(e_1e_2e_3) \qj &= & \qj (e_1e_2e_3)\> =\> -e_2,\\
(e_1e_2e_3) \qk &= & \qk (e_1e_2e_3)\> =\> -e_3,\\
(e_1e_2e_3) e_1 &= & -\qi,\\
(e_1e_2e_3) e_2 &= & -\qj,\\
(e_1e_2e_3) e_3 &= & -\qk.
\end{eqnarray*}
Then,  if $X = a 1 + b \qi + c \qj + d \qk\in \mSpin{3}$,
for every $v = v_1 e_1 + v_2 e_2 + v_3 e_3$, we have
\begin{eqnarray*}
\alpha(X) v X^{-1} & = &
X(v_1 e_1 + v_2 e_2 + v_3 e_3) X^{-1}\\ 
& = & X(e_1e_2e_3)^2(v_1 e_1 + v_2 e_2 + v_3 e_3) X^{-1} \\
& = & (e_1e_2e_3)X (e_1e_2e_3)(v_1 e_1 + v_2 e_2 + v_3 e_3) X^{-1} \\
& = & -(e_1e_2e_3)X (v_1 \qi + v_2 \qj + v_3 \qk) X^{-1}.
\end{eqnarray*}
This shows that the rotation $\rho_X\in \mSO{3}$
induced by $X\in \mSpin {3}$ can be viewed as the rotation
induced by the quaternion 
$a \qone + b \qi + c \qj + d\qk$ on the pure quaternions,
using the maps
\[v \mapsto -(e_1e_2e_3) v, \quad X \mapsto -(e_1e_2e_3) X\]
to go from a vector $v = v_1 e_1 + v_2 e_2 + v_3 e_3$
to the pure quaternion $v_1 \qi + v_2 \qj + v_3 \qk$, and back.

\medskip
We close this section by taking a closer look at
$\mSpin{4}$.
The group  $\mSpin{4}$ consists of all products
\[\prod_{i = 1}^{2k} (a_i e_1 + b_i e_2 + c_i e_3 + d_i e_4)\]
consisting of an even number of factors
and such that $a_i^2 + b_i^2 + c_i^2 + d_i^2  = 1$.
Using the relations
\[e_i^2 = -1,\quad e_je_j = -e_je_i,\quad 1\leq i, j\leq 4,\> i\not= j,\]
every element of $\mSpin{4}$ can be written as
\[x = a_1 1 + a_2 e_1e_2 +  a_3 e_2e_3 +  a_4 e_3e_1 + 
 a_5 e_4e_3 +  a_6 e_4e_1 +  a_7 e_4e_2 + a_8 e_1e_2e_3e_4,\]
where $x$ satisfies the conditions that $xvx^{-1}\in \reals^4$
for all $v\in \reals^4$, and  $N(x) = 1$.
Let 
\[\qi = e_1e_2,\>
  \qj = e_2e_3,\>
  \qk = 3_3e_1,\>
  \qi' = e_4e_3,\>
  \qj' = e_4e_1,\>
  \qk' = e_4e_2,
\]
and
\def\qI{\mathbb{I}}
$\qI = e_1e_2e_3e_4$.
The reader will easily verify that
\begin{eqnarray*}
\qi\qj & = & \qk\\
\qj\qk & = & \qi\\
\qk\qi & = & \qj\\
\qi^2  & = & -1,\quad \qj^2\> =\> -1,\quad \qk^2\> =\> -1\\
\qi\qI & = & \qI\qi\> = \>\qi'\\
\qj\qI & = & \qI\qj\> = \>\qj'\\
\qk\qI & = & \qI\qk\> = \>\qk'\\
\qI^2  & = & 1,\quad \overline{\qI}\> =\> \qI. 
\end{eqnarray*}
Then, every $x\in \mSpin{4}$ can be written as
\[x = u + \qI v,\quad
\hbox{with}\quad u = a 1 + b\qi + c\qj + d\qk
\quad\hbox{and}\quad
v = a' 1 + b'\qi + c'\qj + d'\qk,\]
with the extra conditions stated above.
Using the above identities, we have
\[(u + \qI v)(u' + \qI v') = uu' + vv' + \qI(uv' + vu').\]
As a consequence,
\[N(u + \qI v) = (u + \qI v)(\overline{u} + \qI\overline{v}) =
u\overline{u} + v\overline{v}
+ \qI(u\overline{v} + v\overline{u}),\]
and thus, $N(u + \qI v) = 1$ is equivalent to
\[u\overline{u} + v\overline{v} = 1
\quad\hbox{and}\quad
u\overline{v} + v\overline{u} = 0.\]
As in the case $n = 3$,
it turns out that the conditions $x\in \mCl^0_{4}$
and $N(x) = 1$ imply
that  $xvx^{-1}\in \reals^4$ for all $v\in \reals^4$. 
The only change to the proof is that  $xvx^{-1}\in \mCl^1_{4}$ 
can be written as
\[xvx^{-1} = u + \sum_{i, j, k}\lambda_{i, j, k} e_i e_j e_k,
\quad\hbox{for some $u\in \reals^4$},\quad\hbox{with}\>
\{i, j, k\}\subseteq \{1, 2, 3, 4\}. 
\]
As in the previous proof, we get $\lambda_{i, j, k} = 0$.
Then,
$\mSpin{4}$ consists of those elements
$u + \qI v$ so that
\[u\overline{u} + v\overline{v} = 1
\quad\hbox{and}\quad
u\overline{v} + v\overline{u} = 0,\]
with $u$ and $v$ of the form $a 1 + b\qi + c\qj + d\qk$.
Finally,  we see that $\mSpin{4}$ is isomorphic to
$\mSpin{3}\times\mSpin{3}$ under the isomorphism
\[u + v\qI \mapsto (u + v, u - v).\]
Indeed, we have
\[N(u + v) = (u + v)(\overline{u} + \overline{v}) = 1,\]
and
\[N(u - v) = (u - v)(\overline{u} - \overline{v}) = 1,\]
since
\[u\overline{u} + v\overline{v} = 1
\quad\hbox{and}\quad
u\overline{v} + v\overline{u} = 0,\]
and
\[(u + v, u - v)(u' + v', u' - v') = 
(uu' + vv' + uv' + vu', uu' + vv' - (uv' + vu')).\] 

\remark
It can be shown that the assertion
if $x\in \mCl_{n}^0$ and $N(x) = 1$, then
$xvx^{-1} \in \reals^n$ for all $v\in \reals^n$,
is true up to $n = 5$
(see Porteous \cite{Porteous}, Chapter 13, Proposition 13.58).
However, this is already false for $n = 6$.
For example, if $X = 1/\sqrt{2}(1 + e_1e_2e_3e_4e_5e_6)$,
it is easy to see that $N(X) = 1$, and yet,
$Xe_1X^{-1} \notin \reals^6$.

\section{The Groups $\mPin{p, q}$ and $\mSpin{p, q}$}
\label{seccliff4}
For every nondegenerate quadratic form $\Phi$
over $\reals$,  there is an orthogonal basis
with respect to which $\Phi$ is given by
\[\Phi(x_1, \ldots, x_{p + q}) =
x_1^2 + \cdots + x_p^2 - (x_{p+1}^2 + \cdots + x_{p+q}^2),\]
where $p$ and $q$ only depend on $\Phi$. 
The quadratic form corresponding to $(p, q)$ is denoted
$\Phi_{p, q}$ and we call  $(p, q)$
the  {\it signature of $\Phi_{p, q}$\/}. Let $n = p + q$.
We  define the
group $\mO{p, q}$ as the
group of isometries w.r.t. $\Phi_{p, q}$, {\it i.e.},
the group of linear maps $f$ so that
\[\Phi_{p, q}(f(v)) = \Phi_{p, q}(v)\quad \hbox{for all $v\in \reals^n$}\]
and  the group $\mSO{p, q}$ as the subgroup of $\mO{p, q}$
consisting of the isometries $f\in \mO{p, q}$
with $\det(f) = 1$.
We denote the Clifford algebra $\mCl(\Phi_{p, q})$ where
$\Phi_{p, q}$ has signature $(p, q)$ by
$\mCl_{p, q}$, the corresponding Clifford group
by $\Gamma_{p, q}$, and the special Clifford group
$\Gamma_{p, q}\cap \mCl_{p, q}^0$ by  $\Gamma_{p, q}^+$.
Note that with this new notation,
$\mCl_n = \mCl_{0, n}$.

\danger
As we mentioned earlier, since
Lawson and Michelsohn \cite{Lawson} adopt 
the opposite of our sign convention in defining 
Clifford algebras;
their $\mCl(p, q)$ is our $\mCl(q, p)$.

\medskip
As we mentioned in Section \ref{seccliff2},
we have the problem that $N(v) = -\Phi(v)\cdot 1$,
but $-\Phi(v)$ is not necessarily positive
(where $v\in \reals^n$).
The fix is simple: Allow elements $x\in \Gamma_{p, q}$
with $N(x) = \pm 1$.

\begin{definition}
\label{spindefin2}
We define the {\it pinor group $\mPin{p, q}$\/} as the 
group
\[\mPin{p, q} = \{x\in \Gamma_{p, q} \mid N(x) = \pm 1\},\]
and the 
{\it spinor group $\mSpin{p, q}$\/}
as $\mPin{p, q}\cap \Gamma^+_{p, q}$.
\end{definition}

\remarks
\begin{enumerate}
\item[(1)]
It is easily checked that the group $\mSpin{p, q}$
is also given by
\[\mSpin{p, q} = \{x\in \mCl^0_{p, q} \mid
x v \overline{x} \in \reals^n\quad\hbox{for all $v\in \reals^n$,}\quad
N(x) = 1\}.\]
This is because  $\mSpin{p, q}$ consists of elements of
even degree. 
\item[(2)]
One can check that if $N(x) \not= 0$, then
\[\alpha(x) v x^{-1}= x v t(x)/ N(x).\]
Thus, we have
\[\mPin{p, q} = \{x\in \mCl_{p, q} \mid
x v t(x) N(x) \in \reals^n\quad\hbox{for all $v\in \reals^n$,}\quad
N(x) = \pm 1\}.\]
When $\Phi(x) = -\norme{x}^2$, we have $N(x) = \norme{x}^2$, and
\[\mPin{n} = \{x\in \mCl_{n} \mid
x v t(x) \in \reals^n\quad\hbox{for all $v\in \reals^n$,}\quad
N(x) = 1\}.\]
\end{enumerate}

\medskip
Theorem \ref{spinthm} generalizes as follows:

\begin{theorem}
\label{spinthm2}
The restriction of $\mapdef{\rho}{\Gamma_{p,q}}{\mathbf{GL}(n)}$ 
to the pinor group $\mPin{p, q}$ is a
surjective homomorphism
$\rho\co \mPin{p, q}$ $\rightarrow \mO{p, q}$
whose kernel is $\{-1, 1\}$, and 
the restriction of $\rho$ to the spinor group
$\mSpin{p, q}$ is a
surjective homomorphism
$\mapdef{\rho}{\mSpin{p, q}}{\mSO{p, q}}$
whose kernel is $\{-1, 1\}$.
\end{theorem}

\begin{proof}
The Cartan-Dieudonn\'e also holds for any
nondegenerate quadratic form $\Phi$,
in the sense that every isometry
in $\mO{\Phi}$ is the composition of
reflections defined by hyperplanes orthogonal
to non-isotropic vectors  
(see Dieudonn\'e \cite{Dieudonne67},
Chevalley \cite{Chevalleyv2}, Bourbaki \cite{BourbakiA3}, or 
Gallier \cite{Gallbook2}, Chapter 7, Problem 7.14).
Thus, Theorem \ref{spinthm} also holds
for any nondegenerate quadratic form $\Phi$.
The only change to the proof is the following:
Since $N(w_j) = -\Phi(w_j) \cdot 1$,
we can replace $w_j$ by $w_j/\sqrt{|\Phi(w_j)|}$, so that
$N(w_1\cdots w_k) = \pm 1$, and then
\[f = \rho(w_1\cdots w_k),\]
and $\rho$ is surjective.
\end{proof}

\medskip
If we consider $\reals^n$ equipped with the quadratic form
$\Phi_{p, q}$ (with $n = p + q$), we denote
the set of elements $v\in \reals^n$ with $N(v) = 1$
by $S^{n-1}_{p, q}$.
We have the following corollary of 
Theorem \ref{spinthm2} (generalizing Corollary 
\ref{spincor2}):

\begin{cor}
\label{spincor2}
The group $\mPin{p, q}$ is generated by $S^{n-1}_{p, q}$, 
and every element of $\mSpin{p, q}$ can be written
as the product of an even number of elements
of $S^{n - 1}_{p, q}$.
\end{cor}

\begin{example}
\label{ex3}
The reader should check that
\[\mCl_{0, 1} \approx \complex,\quad \mCl_{1, 0} \approx \reals\oplus \reals.\]
We also have
\[\mPin{0, 1} \approx \integs/4\integs,\quad 
\mPin{1, 0} \approx \integs/2\integs\times \integs/2\integs,\]
from which we get
$\mSpin{0 ,1} = \mSpin{1, 0} \approx \integs/2\integs$.
Also, show that
\[\mCl_{0, 2} \approx \quat,\quad \mCl_{1, 1} \approx M_2(\reals),
\quad \mCl_{2, 0} \approx M_2(\reals),\]
where $M_n(\reals)$ denotes the algebra of $n\times n$ matrices.
One can also work out what are $\mPin{2, 0}$,
$\mPin{1, 1}$, and $\mPin{0, 2}$;
see Choquet-Bruhat \cite{ChoquetBru2}, Chapter I, Section 7, page 26.
Show that
\[\mSpin{0, 2} = \mSpin{2,0} \approx \mU{1},\]
and
\[\mSpin{1, 1} = \{a 1 + b e_1e_2\mid a^2 - b^2 = 1\}.\]
Observe that $\mSpin{1, 1}$ is not connected.
\end{example}

\medskip
More generally, it can be shown that
$\mCl_{p, q}^0$ and $\mCl_{q, p}^0$ are isomorphic, from which
it follows that $\mSpin{p, q}$ and $\mSpin{q, p}$
are isomorphic, but 
$\mPin{p, q}$ and $\mPin{q, p}$ are not isomorphic in general,
and in particular, $\mPin{p, 0}$ and $\mPin{0, p}$ are not isomorphic in general
(see Choquet-Bruhat \cite{ChoquetBru2}, Chapter I, Section 7).
However, due to the ``$8$-periodicity'' of the Clifford algebras
(to be discussed in the next section), it follows that
$\mCl_{p, q}$ and $\mCl_{q, p}$ are isomorphic when
$|p - q| = 0\> \mod\> 4$.

\section{Periodicity of the Clifford Algebras $\mCl_{p, q}$}
\label{seccliff5}
It turns out that the real algebras $\mCl_{p, q}$ can be build up
as tensor products of the basic algebras $\reals$,
$\complex$, and $\quat$. As pointed out by Lounesto
(Section 23.16 \cite{Lounesto}), the description of the 
real algebras $\mCl_{p, q}$ as matrix algebras
and the $8$-periodicity was first observed by Elie Cartan in 1908;
see Cartan's article, {\it Nombres Complexes}, 
based on the original article in German by E. Study,
in Molk \cite{MolkI1},  article I-5 (fasc. 3), pages 329-468.
These algebras are defined in Section 36 under the name
```Systems of Clifford and Lipschitz,'' page 463-466.  
Of course, Cartan used a very different notation;
see page 464 in the article cited above.
These facts were rediscovered independently by Raoul Bott
in the 1960's (see Raoul Bott's comments in Volume 2 of his
Collected papers.).

\medskip
We will use the notation
$\reals(n)$ (resp. $\complex(n)$) for the algebra $M_n(\reals)$ of
all $n \times n$ real matrices  
(resp.  the algebra $M_n(\complex)$ of
all $n \times n$ complex matrices).
As mentioned in Example \ref{ex3}, it is not hard to show that
\begin{eqnarray*}
\mCl_{0, 1} & = & \complex\qquad \mCl_{1, 0} = \reals\oplus \reals\\
\mCl_{0, 2} & = & \quat\qquad \mCl_{2, 0} = \reals(2),
\end{eqnarray*}
and
\[\mCl_{1, 1} = \reals(2).\]
The key to the classification is the following
lemma:
\begin{lemma}
\label{periodlem}
We have the isomorphisms
\begin{eqnarray*}
\mCl_{0, n+2} &\approx & \mCl_{n, 0}\tensor \mCl_{0, 2} \\
\mCl_{n+2, 0} &\approx & \mCl_{0, n}\tensor \mCl_{2, 0} \\
\mCl_{p+1, q+1} &\approx & \mCl_{p, q}\tensor \mCl_{1, 1},
\end{eqnarray*}
for all $n, p, q\geq 0$.
\end{lemma}

\begin{proof}
Let $\Phi_{0, n}(x) = -\norme{x}^2$, where $\norme{x}$ is the
standard Euclidean norm on $\reals^{n+2}$, and
let $(e_1, \ldots, e_{n+2})$ be an orthonormal basis for 
$\reals^{n+2}$ under the standard 
Euclidean inner product. 
We also let $(e_1', \ldots, e_{n}')$ be a set
of generators for $\mCl_{n, 0}$ and
$(e_1'', e_2'')$ be a set
of generators for $\mCl_{0, 2}$.
We can define a linear map 
$\mapdef{f}{\reals^{n+2}}{\mCl_{n, 0}\tensor \mCl_{0, 2}}$
by its action on the basis
$(e_1, \ldots, e_{n+2})$ as follows:
\[
f(e_i) =
\oldcases{
e_i'\tensor e_1'' e_2'' & for $1 \leq i \leq n$\cr
1 \tensor e_{i - n}'' & for $n+1 \leq i \leq n+2$.\cr
}
\]
Observe that for $1\leq i, j \leq n$, we have
\[f(e_i)f(e_j) + f(e_j)f(e_i) = (e_i'e_j' + e_j'e_i')\tensor (e_1''e_e'')^2
= -2\delta_{i j} 1\tensor 1,\]
since $e_1''e_2'' = -e_2''e_1''$, $(e_1'')^2 = -1$,
and $(e_2'')^2 = -1$, and
$e_i'e_j' = -e_j'e_i'$, for all $i\not= j$, and
$(e_i')^2 = 1$, for all $i$ with $1 \leq i \leq n$.
Also, for $n + 1 \leq i, j \leq n + 2$, we have
\[f(e_i)f(e_j) + f(e_j)f(e_i) 
= 1\tensor (e_{i-n}''e_{j-n}'' + e_{j-n}''e_{i-n}'')
= -2\delta_{i j} 1\tensor 1,\]
and
\[f(e_i)f(e_k) + f(e_k)f(e_i) = 2 e_i'\tensor
(e_1''e_2''e_{n - k}'' + e_{n - k}''e_1''e_2'')
= 0,\]
for   $1\leq i, j \leq n$ and  $n + 1 \leq k \leq n + 2$
(since $e_{n - k}'' = e_1''$ or $e_{n - k}'' = e_2''$).
Thus, we have
\[f(x)^2 = -\norme{x}^2\cdot 1\tensor 1\quad
\hbox{for all $x\in \reals^{n+2}$},
\] 
and by the universal mapping property
of $\mCl_{0, n+2}$, we get an algebra map
\[\mapdef{\widetilde{f}}{\mCl_{0, n+2}}{\mCl_{n, 0}\tensor \mCl_{0, 2}}.\]
Since $\widetilde{f}$ maps onto a set of
generators, it is surjective. However
\[\dimm(\mCl_{0, n+2}) = 2^{n + 2} = 2^n \cdot 2 =
\dimm(\mCl_{n, 0})\dimm(\mCl_{0, 2}) =
\dimm(\mCl_{n, 0}\tensor \mCl_{0, 2}),\]
and $\widetilde{f}$ is an isomorphism.

\medskip
The proof of the second identity is analogous.
For the third identity, we have
\[\Phi_{p, q}(x_1, \ldots, x_{p + q}) =
x_1^2 + \cdots + x_p^2 - (x_{p+1}^2 + \cdots + x_{p+q}^2),\]
and let 
$(e_1, \ldots, e_{p+1}$, $\epsilon_1, \ldots, \epsilon_{q+1})$ 
be an orthogonal basis for $\reals^{p + q + 2}$ so that
$\Phi_{p+1, q+1}(e_i) = +1$ and $\Phi_{p+1, q+1}(\epsilon_j) = -1$ 
for $i = 1, \ldots, p+1$ and $j = 1, \ldots, q+1$.
Also, let 
$(e_1', \ldots, e_{p}'$, $\epsilon_1', \ldots, \epsilon_{q}')$
be a set of generators for $\mCl_{p, q}$
and $(e_1'', \epsilon_1'')$ be a set of generators
for $\mCl_{1, 1}$.
We define a linear map
$\mapdef{f}{\reals^{p+q+2}}{\mCl_{p, q}\tensor \mCl_{1, 1}}$
by its action on the basis
as follows:
\[
f(e_i) =
\oldcases{
e_i'\tensor e_1'' \epsilon _1'' & for $1 \leq i \leq p$\cr
1 \tensor e_1'' & for $i = p+1$,\cr
}
\]
and
\[
f(\epsilon_j) =
\oldcases{
\epsilon_j'\tensor e_1'' \epsilon _1'' & for $1 \leq j \leq q$\cr
1 \tensor \epsilon_1'' & for $j = q+1$.\cr
}
\]
We can check that
\[f(x)^2 = \Phi_{p+1, q+1}(x)\cdot 1\tensor 1\quad
\hbox{for all $x\in \reals^{p + q + 2}$},\]
and we finish the proof as in the first case.
\end{proof}

\medskip
To apply this lemma, we need some further isomorphisms
among various matrix algebras.

\begin{proposition}
\label{matiso}
The following isomorphisms hold:
\begin{eqnarray*}
\reals(m)\tensor \reals(n) &\approx & \reals(mn)
\qquad \hbox{for all $m, n\geq 0$}\\
\reals(n)\tensor_{\reals} K & \approx & K(n)
\qquad \hbox{for $K = \complex$ or $K = \quat$ and all $n\geq 0$}\\
\complex\tensor_{\reals} \complex & \approx & \complex\oplus \complex\\
\complex\tensor_{\reals} \quat & \approx & \complex(2)\\
\quat\tensor_{\reals} \quat & \approx & \reals(4).
\end{eqnarray*}
\end{proposition}

\begin{proof}
Details can be found in Lawson and Michelsohn \cite{Lawson}.
The first two isomorphisms are quite obvious. The third isomorphism
$\complex \oplus \complex \rightarrow \complex\tensor \complex$
is obtained by sending
\[(1, 0) \mapsto \frac{1}{2}(1\tensor 1 + i\tensor i),\quad
(0, 1) \mapsto \frac{1}{2}(1\tensor 1 - i\tensor i).\]
The field $\complex$ is isomorphic to the
subring of $\quat$ generated by $\qi$. Thus,
we can view $\quat$ as a $\complex$-vector space
under left scalar multiplication.
Consider the $\reals$-bilinear map \\
$\mapdef{\pi}{\complex\times \quat}{\Hom{\complex}{\quat}{\quat}}$
given by
\[\pi_{y,z}(x) = yx\overline{z},\]
where $y\in \complex$ and $x, z\in \quat$.
Thus, we get an $\reals$-linear map  
$\mapdef{\pi}{\complex\tensor_{\reals}\quat}{\Hom{\complex}{\quat}{\quat}}$.
However, we have
$\Hom{\complex}{\quat}{\quat}\approx \complex(2)$.
Furthermore, since
\[\pi_{y, z}\circ \pi_{y', z'} = \pi_{yy', zz'},\]
the map $\pi$ is an algebra homomorphism
\[\mapdef{\pi}{\complex\times \quat}{\complex(2)}.\]
We can check on a basis that $\pi$ is injective,
and since
\[\dimm_{\reals}(\complex\times \quat) = \dimm_{\reals}(\complex(2)) = 8,\]
the map $\pi$ is an isomorphism.
The last isomorphism is proved in a similar fashion.
\end{proof}

\medskip
We now have the main periodicity theorem.

\begin{theorem}
\label{periodthm} (Cartan/Bott)
For all $n\geq 0$, we have the following isomorphisms:
\begin{eqnarray*}
\mCl_{0, n+8} & \approx & \mCl_{0, n}\tensor \mCl_{0, 8}\\
\mCl_{n+8, 0} & \approx & \mCl_{n, 0}\tensor \mCl_{8, 0}.
\end{eqnarray*}
Furthermore,
\[\mCl_{0, 8} = \mCl_{8, 0} = \reals(16).\]
\end{theorem}

\begin{proof}
By Lemma \ref{periodlem}
we have the isomorphisms
\begin{eqnarray*}
\mCl_{0, n+2} &\approx & \mCl_{n, 0}\tensor \mCl_{0, 2} \\
\mCl_{n+2, 0} &\approx & \mCl_{0, n}\tensor \mCl_{2, 0}, 
\end{eqnarray*}
and thus, 
\[\mCl_{0, n+8} \approx  \mCl_{n+6, 0}\tensor \mCl_{0, 2} 
\approx \mCl_{0, n+4}\tensor \mCl_{2, 0}  \tensor \mCl_{0, 2} \approx\cdots
\approx \mCl_{0, n}\tensor \mCl_{2, 0}  \tensor \mCl_{0, 2} \tensor
\mCl_{2, 0}  \tensor \mCl_{0, 2}.\]
Since $\mCl_{0,2} = \quat$ and $\mCl_{2, 0} = \reals(2)$,
by Proposition \ref{matiso}, we get
\[\mCl_{2, 0}  \tensor \mCl_{0, 2} \tensor
\mCl_{2, 0}  \tensor \mCl_{0, 2} \approx \quat\tensor\quat\tensor
\reals(2)\tensor\reals(2) \approx \reals(4)\tensor \reals(4) \approx
\reals(16).\]
The second isomorphism is proved in a similar fashion.
\end{proof}

\medskip
From all this, we can deduce the following table:
\[
\begin{matrix}
 n          &     0  & 1  &  2  &  3  & 4  & 5 & 6 &  7  &  8 \\
\mCl_{0, n} & \reals & \complex &\quat &\quat\oplus\quat
&\quat(2) & \complex(4) & \reals(8) &\reals(8)\oplus \reals(8) &\reals(16)\\
\mCl_{n, 0} & \reals & \reals\oplus\reals &\reals(2) &\complex(2)
&\quat(2) & \quat(2)\oplus\quat(2) & \quat(4) &\complex(8) &\reals(16)
\end{matrix}
\]

\medskip
A table of the Clifford groups $\mCl_{p, q}$
for $0 \leq p, q\leq 7$ can be found in Kirillov \cite{Kirillov01},
and for $0 \leq p, q\leq 8$,
in Lawson and Michelsohn \cite{Lawson}
(but beware that their $\mCl_{p, q}$ is our  $\mCl_{q, p}$).
It can also be shown that
\begin{eqnarray*}
\mCl_{p+1, q} & \approx & \mCl_{q+1, p} \\
\mCl_{p, q} & \approx & \mCl_{p - 4, q + 4}
\end{eqnarray*}
with $p\geq 4$ in the second identity
(see Lounesto  \cite{Lounesto}, Chapter 16, Sections 16.3 and 16.4).
Using the second identity, if $|p - q| = 4k$,
it is easily shown by induction on $k$ that
$\mCl_{p,q} \approx \mCl_{q, p}$, as claimed in the previous
section.

\medskip
We also have the isomorphisms
\[\mCl_{p, q} \approx \mCl^0_{p, q+1},\]
frow which it follows that
\[\mSpin{p, q} \approx \mSpin{q, p}\]
(see Choquet-Bruhat \cite{ChoquetBru2}, Chapter I, Sections 4 and 7).
However, in general, $\mPin{p,q}$ and $\mPin{q, p}$ are not isomorphic.
In fact,  $\mPin{0,n}$ and $\mPin{n, 0}$ are not isomorphic
if $n \not= 4k$, with $k\in \natnums$
(see Choquet-Bruhat \cite{ChoquetBru2}, Chapter I, Section 7, page 27).

\section{The Complex Clifford Algebras $\mCl(n, \complex)$}
\label{seccliff6}
One can also consider Clifford algebras over the complex field $\complex$.
In this case, it is well-known that every nondegenerate quadratic form
can be expressed by
\[\Phi_n^{\complex}(x_1, \ldots, x_n) = x_1^2 + \cdots + x_n^2\]
in some orthonormal basis.
Also, it is easily shown that the complexification
$\complex\tensor_{\reals}\mCl_{p, q}$ of the real Clifford
algebra $\mCl_{p, q}$ is isomorphic to
$\mCl(\Phi^{\complex}_{n})$. Thus, all these complex algebras
are isomorphic for $p + q = n$, and we denote them by $\mCl(n, \complex)$.
Theorem \ref{periodlem} yields the following periodicity theorem:

\begin{theorem}
\label{periodthmc}
The following isomorphisms hold:
\[\mCl(n+2, \complex) \approx \mCl(n, \complex) \tensor_{\complex} \mCl(2, \complex),\]
with $\mCl(2, \complex) = \complex(2)$.
\end{theorem}

\begin{proof}
Since $\mCl(n, \complex) = \complex\tensor_{\reals} \mCl_{0, n}
=  \complex\tensor_{\reals}\mCl_{n, 0}$,
by Lemma \ref{periodlem}, we have
\[\mCl(n+2, \complex) = \complex\tensor_{\reals} \mCl_{0, n+2}
 \approx \complex\tensor_{\reals}(\mCl_{n, 0}\tensor_{\reals} \mCl_{0, 2})
\approx  (\complex\tensor_{\reals}\mCl_{n, 0}) \tensor_{\complex}
(\complex\tensor_{\reals} \mCl_{0, 2}) 
.\]
However, $\mCl_{0, 2} = \quat$, 
 $\mCl(n, \complex) = \complex\tensor_{\reals}\mCl_{n, 0}$,
and
$\complex \tensor_{\reals}\quat \approx \complex(2)$, so we get
$\mCl(2, \complex) = \complex(2)$ and
\[\mCl(n+2, \complex) \approx \mCl(n, \complex) \tensor_{\complex} \complex(2),\]
and the theorem is proved.
\end{proof}

\medskip
As a corollary of Theorem \ref{periodthmc}, we obtain the fact
that
\[\mCl(2k, \complex) \approx \complex(2^k)
\quad\hbox{and}\quad
\mCl(2k+1, \complex) \approx \complex(2^k)\oplus \complex(2^k).\]
The table of the previous section can also be completed as follows:
\[
\begin{matrix}
 n          &     0  & 1  &  2  &  3  & 4  & 5 & 6 &  7  &  8 \\
\mCl_{0, n} & \reals & \complex &\quat &\quat\oplus\quat
&\quat(2) & \complex(4) & \reals(8) &\reals(8)\oplus \reals(8) &\reals(16)\\
\mCl_{n, 0} & \reals & \reals\oplus\reals &\reals(2) &\complex(2)
&\quat(2) & \quat(2)\oplus\quat(2) & \quat(4) &\complex(8) &\reals(16)\\
\mCl(n,\complex) & \complex & 2\complex & \complex(2) & 
2\complex(2) & \complex(4) & 2\complex(4) &
\complex(8) & 2\complex(8) & \complex(16). 
\end{matrix}
\]
where $2\complex(k)$ is an abbrevation for $\complex(k) \oplus \complex(k)$.

\section[The Groups $\mPin{p, q}$ and $\mSpin{p, q}$ as double covers]
{The Groups $\mPin{p, q}$ and $\mSpin{p, q}$ as double covers
of $\mO{p, q}$ and $\mSO{p, q}$}
\label{seccliff7}
It turns out that the groups $\mPin{p, q}$ and $\mSpin{p, q}$ 
have nice topological properties w.r.t. the groups
$\mO{p, q}$ and $\mSO{p, q}$. To explain this, we review
the definition of covering maps and covering spaces
(for details, see Fulton \cite{Fulton95}, Chapter 11).
Another interesting source is Chevalley \cite{Chevalley46},
where is is proved that $\mSpin{n}$ is
a universal double cover of $\mSO{n}$ for all $n\geq 3$.

\medskip
Since $C_{p, q}$ is an algebra of dimension
$2^{p + q}$, it is a topological space as a vector space isomorphic
to $V = \reals^{2^{p+q}}$. Now, the group  $C^*_{p, q}$ of units of $C_{p, q}$ is
open in $C_{p, q}$, because $x\in C_{p, q}$ is a unit if the linear
map $y \mapsto xy$ is an isomorphism,
and $\mGL{V}$ is open in $\mathrm{End}(V)$, the space of endomorphisms
of $V$.
Thus, $C^*_{p, q}$ is a Lie group, and since 
$\mPin{p, q}$ and $\mSpin{p, q}$ are clearly closed subgroups of
$C^*_{p, q}$, they are also Lie groups.


\begin{definition}
\label{covermap}
Given two topological spaces $X$ and $Y$, a {\it covering map\/}
is a continuous surjective map $\mapdef{p}{Y}{X}$ with the property
that for every $x\in X$, there is some open subset $U\subseteq X$
with $x\in U$, so that $p^{-1}(U)$ is the disjoint union of
open subsets $V_{\alpha}\subseteq Y$, and the restriction of $p$
to each $V_{\alpha}$ is a homeomorphism onto $U$. We say that
$U$ is {\it evenly covered by $p$\/}.
We also say that $Y$ is a {\it covering space of $X$\/}.
A covering map  $\mapdef{p}{Y}{X}$ is called {\it trivial\/}
if $X$ itself is evenly covered by $p$ ({\it i.e.},
$Y$ is the disjoint union of open subsets $Y_{\alpha}$
each homeomorphic to $X$), and {\it nontrivial\/} otherwise.
When each fiber $p^{-1}(x)$ has the same finite cardinaly $n$
for all $x\in X$, we say that $p$ is an {\it $n$-covering\/}
(or {\it $n$-sheeted covering\/}).
\end{definition}

\medskip
Note that a covering map  $\mapdef{p}{Y}{X}$ is not always trivial,
but always {\it locally trivial\/} ({\it i.e.},
for every $x\in X$, it is trivial in some open neighborhood of $x$).
A covering is trivial iff $Y$ is isomorphic
to a product space of $X\times T$, where
$T$ is any set with the discrete topology.
Also, if $Y$ is connected, then the covering map is nontrivial.

\begin{definition}
\label{mapcovermap}
An {\it isomorphism\/} $\varphi$
between covering maps $\mapdef{p}{Y}{X}$ and $\mapdef{p'}{Y'}{X}$
is a homeomorphism $\mapdef{\varphi}{Y}{Y'}$ so that
$p = p' \circ \varphi$.
\end{definition}

\medskip
Typically, the space $X$ is connected, in which case it can be shown
that all the fibers $p^{-1}(x)$ have the same cardinality.

\medskip
One of the most important properties of covering spaces
is the path--lifting property, a property that we will use to 
show that  $\mSpin{n}$ is path-connected. 

\begin{proposition}
\label{pathlift} (Path lifting)
Let $\mapdef{p}{Y}{X}$ be a covering map, and let
$\mapdef{\gamma}{[a, b]}{X}$ be any continuous
curve from $x_a = \gamma(a)$ to $x_b = \gamma(b)$
in $X$. If $y\in Y$ is any point so that $p(y) = x_a$, then
there is a unique curve
$\mapdef{\widetilde{\gamma}}{[a, b]}{Y}$ so that 
$y = \widetilde{\gamma}(a)$ and
\[p \circ \widetilde{\gamma}(t) = \gamma(t)\quad
\hbox{for all $t\in [a, b]$}.\]
\end{proposition}

\begin{proof}
See Fulton \cite{Fulton91}, Chapter 11, Lemma 11.6.
\end{proof}

\medskip
Many important covering maps arise from the action of a group
$G$ on a space $Y$.
If $Y$ is a topological space, an {\it action (on the left) of
a group $G$ on $Y$\/} is a map $\mapdef{\alpha}{G\times Y}{Y}$
satisfying the following conditions,
where for simplicity of notation, we denote $\alpha(g, y)$
by $g\cdot y$:
\begin{enumerate}
\item[(1)]
$g\cdot (h \cdot y) = (gh)\cdot y$, for all $g, h\in G$ and $y\in Y$;
\item[(2)]
$1\cdot y = y$, for all $\in Y$, where $1$ is the identity of the group $G$;
\item[(3)]
The map $y \mapsto g\cdot y$ is a homeomorphism of $Y$ for
every $g\in G$.
\end{enumerate}

We define an equivalence relation on $Y$ as follows:
$x\equiv y$ iff $y = g\cdot x$ for some $g\in G$
(check that this is an equivalence relation).
The equivalence class $G\cdot x = \{g\cdot x\mid g\in G\}$
of any $x\in Y$ is called the {\it orbit of $x$\/}.
We obtain the quotient space $Y/G$ and 
the projection  map $\mapdef{p}{Y}{Y/G}$ sending
every $y\in Y$ to its orbit.
The space $Y/G$ is given the quotient topology 
(a subset $U$ of $Y/G$ is open iff $p^{-1}(U)$ is open in $Y$).

\medskip
Given a subset $V$ of $Y$ and any $g\in G$, we let
\[g\cdot V = \{g\cdot y\mid y\in V\}.\]
We say that $G$ {\it acts evenly on $Y$\/} if for every $y\in Y$,
there is an open subset $V$ containing $y$ so that
$g\cdot V$ and $h\cdot V$ are disjoint for any two
distinct elements $g, h\in G$.

\medskip
The importance of the notion a group acting evenly 
is that such actions induce a covering map.

\begin{proposition}
\label{gactcov}
If $G$ is a group acting evenly on a space $Y$, then
the projection map $\mapdef{p}{Y}{Y/G}$ is a covering map.
\end{proposition}

\begin{proof}
See Fulton \cite{Fulton91}, Chapter 11, Lemma 11.17.
\end{proof}

\medskip
The following proposition shows that $\mPin{p, q}$ and $\mSpin{p, q}$ 
are nontrivial covering spaces, unless $p = q = 1$.

\begin{proposition}
\label{doublecov}
For all $p, q\geq 0$, the groups
$\mPin{p, q}$ and $\mSpin{p, q}$  are double covers
of  $\mO{p, q}$ and $\mSO{p, q}$, respectively.
Furthermore, they are nontrivial covers
unless $p = q = 1$.
\end{proposition}

\begin{proof}
We know that kernel of the homomorphism
$\mapdef{\rho}{\mPin{p, q}}{\mO{p, q}}$ is
$\integs_2 = \{-1, 1\}$. If we let $\integs_2$ act
on $\mPin{p, q}$ in the natural way, then
$\mO{p, q} \approx \mPin{p, q}/\integs_2$, and
the reader can easily check that $\integs_2$ acts
evenly. 
By Proposition \ref{gactcov}, we get a double cover.
The argument for $\mapdef{\rho}{\mSpin{p, q}}{\mSO{p, q}}$ is
similar.

\medskip
Let us now assume that $p\not= 1$ or $q\not= 1$.
In order to prove that we have nontrivial covers,
it is enough to show that $-1$ and $1$ are connected
by a path in $\mPin{p, q}$
(If we had $\mPin{p,q} = U_1\cup U_2$
with $U_1$ and $U_2$ open, disjoint, and homeomorphic
to $\mO{p, q}$, then $-1$ and $1$ would not be in the same
$U_i$, and so, they would be in disjoint
connected components. Thus, 
$-1$ and $1$ can't be path--connected,
and similarly with  $\mSpin{p, q}$ and $\mSO{p, q}$.)
Since $(p, q)\not= (1, 1)$, we can find
two orthogonal vectors $e_1$ and $e_2$ so that
$\Phi_{p, q}(e_1) = \Phi_{p, q}(e_2) = \pm 1$.
Then,
\[\gamma(t) = \pm\cos(2t)\, 1 + \sin(2t)\, e_1e_2 =
(\cos t\, e_1  + \sin t\, e_2)(\sin t\, e_2  - \cos t\, e_1),\]
for $0 \leq t\leq \pi$, defines a path in
$\mSpin{p, q}$, since 
\[(\pm\cos(2t)\, 1 + \sin(2t)\, e_1e_2)^{-1}
= \pm\cos(2t)\, 1 - \sin(2t)\, e_1e_2,\] 
as desired.
\end{proof}

\medskip
In particular, if $n \geq 2$, 
since the group $\mSO{n}$ is path-connected,
the group $\mSpin{n}$  is also path-connected.
Indeed, given any two points $x_a$ and  $x_b$ in 
$\mSpin{n}$, there is a path $\gamma$
from $\rho(x_a)$ to $\rho(x_b)$ in $\mSO{n}$ 
(where $\mapdef{\rho}{\mSpin{n}}{\mSO{n}}$
is the covering map). By Proposition \ref{pathlift},
there is a path $\widetilde{\gamma}$ in $\mSpin{n}$ 
with origin $x_a$ and some origin $\widetilde{x_b}$
so that $\rho(\widetilde{x_b}) = \rho(x_b)$.
However, $\rho^{-1}(\rho(x_b)) = \{-x_b, x_b\}$,
and so 
$\widetilde{x_b} = \pm x_b$. The argument
used in the proof of Proposition \ref{doublecov}
shows that $x_b$ and $-x_b$ are path-connected,
and so, there is a path from $x_a$ to $x_b$, and
$\mSpin{n}$ is path-connected.

\medskip
In fact, for $n\geq 3$, it turns out that 
$\mSpin{n}$ is simply connected. Such a
covering space is called a {\it universal cover\/}
(for instance, see Chevalley \cite{Chevalley46}).

\medskip
This last fact requires more algebraic topology than we are 
willing to explain in detail, and we only sketch the proof.
The notions of fibre bundle, fibration, and
homotopy sequence associated with a fibration are
needed in the proof. We refer the perseverant
readers to Bott and Tu \cite{BottTu} 
(Chapter 1 and Chapter 3, Sections 16--17)
or Rotman \cite{Rotman} (Chapter 11)
for a detailed explanation of these concepts.

\medskip
Recall that a topological space is {\it simply connected\/}
if it is path connected and if $\pi_1(X) = (0)$,
which means that every closed path in $X$ is homotopic to a point.
Since we just proved that $\mSpin{n}$ is path connected
for $n \geq 2$, we just need to prove that
$\pi_1(\mSpin{n}) = (0)$ for all $n \geq 3$. 
The following facts are needed to prove the above assertion:
\begin{enumerate}
\item[(1)]
The sphere $S^{n-1}$ is simply connected for all $n \geq 3$.
\item[(2)]
The group $\mSpin{3} \simeq \mSU{2}$ is homeomorphic to 
$S^3$, and thus, $\mSpin{3}$ is simply connected.
\item[(3)]
The group $\mSpin{n}$ acts on $S^{n - 1}$ in such a way that
we have a fibre bundle with fibre $\mSpin{n-1}$:
\[\mSpin{n-1} \longrightarrow \mSpin{n} \longrightarrow S^{n-1}.\]
\end{enumerate}

Fact (1) is a standard proposition of algebraic topology,
and a proof can found in many books.
A particularly elegant and yet simple argument consists
in showing that any closed curve on $S^{n-1}$
is homotopic to a curve that omits some point.
First, it is easy to see that in $\reals^n$, any closed
curve is homotopic to a piecewise linear curve (a polygonal curve), and
the radial projection of such a curve on $S^{n-1}$
provides the desired curve. Then, we use the stereographic
projection  of $S^{n-1}$ from any point omitted by that curve
to get another closed curve in
$\reals^{n-1}$. Since $\reals^{n-1}$ is simply connected,
that curve is homotopic to a point, and so
is its preimage curve on $S^{n-1}$.
Another simple proof uses a special version of the
Seifert---van Kampen's theorem (see Gramain \cite{Gramain}).

\medskip
Fact (2) is easy to establish directly, using (1).

\medskip
To prove (3), we let $\mSpin{n}$ act on $S^{n-1}$
via the standard action: $x \cdot v = xvx^{-1}$.
Because $\mSO{n}$ acts transitively on $S^{n-1}$
and there is a surjection $\mSpin{n} \longrightarrow \mSO{n}$,
the group $\mSpin{n}$ also acts transitively on $S^{n-1}$.
Now, we have to show that the stabilizer of any element of
$S^{n-1}$ is $\mSpin{n-1}$. For example, we can do this for
$e_1$. This amounts to some simple calculations
taking into account the identities among basis elements.
Details of this proof can be found in 
Mneimn\'e  and Testard \cite{Mneimne}, Chapter 4.
It is still necessary to prove that $\mSpin{n}$ is a fibre bundle
over $S^{n-1}$ with fibre $\mSpin{n-1}$. 
For this, we use the following
results whose proof can be found in 
Mneimn\'e  and Testard \cite{Mneimne}, Chapter 4:

\begin{lemma}
\label{fibrelem1}
Given any topological group $G$, if $H$ is a closed
subgroup of $G$ and the projection
$\mapdef{\pi}{G}{G/H}$ has a local section
at every point of $G/H$, then 
\[H \longrightarrow G \longrightarrow G/H\]
is a fibre bundle with fibre $H$.
\end{lemma}

Lemma \ref{fibrelem1} 
implies the following key proposition:

\begin{proposition}
\label{fibreprop2}
Given any linear Lie group $G$, if $H$ is a closed
subgroup of $G$, then 
\[H \longrightarrow G \longrightarrow G/H\]
is a fibre bundle with fibre $H$.
\end{proposition}

\medskip
Now, a fibre bundle is a fibration (as defined in Bott and Tu
\cite{BottTu},
Chapter 3, Section 16, or in Rotman \cite{Rotman}, Chapter 11). 
For a proof of this fact, see Rotman \cite{Rotman}, Chapter 11, or
Mneimn\'e  and Testard \cite{Mneimne}, Chapter 4.
So, there is a homotopy sequence associated with the fibration
(Bott and Tu  \cite{BottTu}, Chapter 3, Section 17, or
Rotman \cite{Rotman}, Chapter 11, Theorem 11.48),
and in particular, we have the exact sequence
\[\pi_1(\mSpin{n-1}) \longrightarrow \pi_1(\mSpin{n}) 
\longrightarrow \pi_1(S^{n-1}).\]
Since $\pi_1(S^{n-1}) = (0)$ for $n \geq 3$, we get a surjection
\[\pi_1(\mSpin{n-1}) \longrightarrow \pi_1(\mSpin{n}),\]
and so, by induction and (2), we get
\[\pi_1(\mSpin{n}) \approx \pi_1(\mSpin{3}) = (0),\]
proving that $\mSpin{n}$ is simply connected for $n\geq 3$.

\medskip
We can also show that $\pi_1(\mSO{n}) = \integs/2\integs$
for all $n \geq 3$. For this, we use Theorem \ref{spinthm}
and Proposition \ref{doublecov}, which imply that
$\mSpin{n}$ is a fibre bundle over $\mSO{n}$ with fibre
$\{-1, 1\}$, for $n \geq 2$:
\[\{-1, 1\} \longrightarrow \mSpin{n} \longrightarrow \mSO{n}.\]
Again, the homotopy sequence of the fibration exists,
and in particular we get the exact sequence
\[\pi_1(\mSpin{n}) \longrightarrow \pi_1(\mSO{n})
\longrightarrow \pi_0(\{-1, +1\}) \longrightarrow \pi_0(\mSO{n}).\]
Since $\pi_0(\{-1, +1\}) = \integs/2\integs$,
$\pi_0(\mSO{n})= (0)$,  and
$\pi_1(\mSpin{n}) = (0)$, when $n \geq 3$,  we get
the exact sequence
\[(0) \longrightarrow \pi_1(\mSO{n})
\longrightarrow \integs/2\integs \longrightarrow (0),\]
and so,  $\pi_1(\mSO{n}) = \integs/2\integs$.
Therefore, $\mSO{n}$ is not simply connected for $n \geq 3$.

\remark
Of course, we have been rather cavalier in our presentation.
Given a topological space $X$, the group
$\pi_1(X)$ is the {\it fundamental group of $X$\/}, {\it i.e.}
the group of homotopy classes of closed paths in $X$
(under composition of loops).
But $\pi_0(X)$ is generally {\it not\/} a group!
Instead, $\pi_0(X)$ is the set of path-connected
components of $X$. However, when $X$ is a Lie group,
$\pi_0(X)$ is indeed a group. Also, we have to make
sense of what it means for the sequence
to be exact.  All this can be made rigorous
(see Bott and Tu  \cite{BottTu}, Chapter 3, Section 17,
or Rotman \cite{Rotman}, Chapter 11).

\section[More on the Topology of $\mO{p, q}$ and $\mSO{p, q}$]
{More on the Topology of $\mO{p, q}$ and $\mSO{p, q}$}
\label{seccliff8}
It turns out that the topology
of the group $\mO{p, q}$ is completely determined by the topology
of $\mO{p}$ and $\mO{q}$. This result can be obtained
as a simple consequence of some standard Lie group theory.
The key notion is that of a pseudo-algebraic group.

\medskip
Consider the group $\mGL{n, \complex}$ of invertible
$n\times n$ matrices with complex coefficients.
If $A = (a_{k l})$ is such a matrix, denote
by $x_{k l}$ the real part (resp. $y_{k l}$, the imaginary
part) of $a_{k l}$ (so, $a_{k l} = x_{k l} + i y_{k l}$).

\begin{defin}
\label{pseudoalgdef}
A subgroup $G$ of $\mGL{n, \complex}$ is {\it pseudo-algebraic\/}
iff there is a finite set of polynomials in $2n^2$ variables
with real coefficients 
$\{P_i(X_1, \ldots, X_{n^2}, Y_1, \ldots, Y_{n^2})\}_{i = 1}^{t}$,
so that 
\[
A = (x_{k l} + i y_{k l})\in G
\quad\hbox{iff}\quad 
P_i(x_{1 1}, \ldots, x_{n n}, y_{1 1}, \ldots, y_{n n}) = 0,
\quad\hbox{for $i = 1,\ldots, t$}.
\]
\end{defin}

\medskip
Recall that if $A$ is a complex $n\times n$-matrix,
its {\it adjoint\/} $A^*$ is defined by
$A^* = \transpos{(\overline{A})}$.
Also, $\mU{n}$ denotes the group of unitary
matrices, {\it i.e.}, those matrices $A\in \mGL{n,\complex}$
so that $AA^* = A^*A = I$, and
$\mH{n}$ denotes the vector space of Hermitian matrices,
{\it i.e.},  those matrices $A$
so that $A^* = A$.
Then, we have the following theorem which is essentially
a refined version of the polar decomposition of
matrices:

\begin{thm}
\label{genpolar}
Let $G$ be a pseudo-algebraic subgroup of $\mGL{n, \complex}$ stable under
adjunction ({\it i.e.}, we have $A^*\in G$ whenever $A\in G$). Then,
there is some integer $d\in \natnums$ so that
$G$ is homeomorphic to $(G \cap \mU{n})\times \reals^d$.
Moreover, if $\mfrac{g}$ is the Lie algebra of $G$, 
the map 
\[(\mU{n}\cap G)\times (\mH{n}\cap \mfrac{g}) 
\longrightarrow G,\quad
\hbox{given by}\quad
(U, H) \mapsto U e^H,
\]
is a homeomorphism onto $G$.
\end{thm}

\begin{proof}
A proof can be found in Knapp \cite{Knapp96}, Chapter 1, or
Mneimn\'e  and Testard \cite{Mneimne}, Chapter 3.
\end{proof}

\medskip
We now apply Theorem \ref{genpolar} to determine the
structure of the space $\mO{p, q}$.
Let $J_{p, q}$ be the matrix
\[
J_{p, q} = \amsmata{I_p}{0}{0}{-I_q}.
\]
We know that $\mO{p, q}$ consists of the matrices
$A$ in $\mGL{p+q, \reals}$ such that
\[
\transpos{A} J_{p, q} A = J_{p, q},
\]
and so $\mO{p, q}$ is clearly pseudo-algebraic.
Using the above equation, it is easy to determine the Lie
algebra, $\mfrac{o}(p, q)$, of $\mO{p, q}$.
We find that $\mfrac{o}(p, q)$
is given by
\[
\mfrac{o}(p, q) = \left\{
\left.
\amsmata{X_1}{X_2}{\transpos{X_2}}{X_3}\> 
\right|
\>
\transpos{X_1} = - X_1,\>\> \transpos{X_3} = - X_3,\>\>
X_2\>\> \hbox{arbitrary} 
\right\},
\]
where $X_1$ is a $p\times p$ matrix,  $X_3$ is a $q\times q$ matrix,
and $X_2$ is a $p\times q$ matrix.
Consequently, it immediately follows that
\[
\mfrac{o}(p, q)\cap \mH{p+q} = \left\{
\left.
\amsmata{0}{X_2}{\transpos{X_2}}{0}\> 
\right|
\>
X_2\>\> \hbox{arbitrary} 
\right\},
\]
a vector space of dimension $pq$.

\medskip
Some simple calculations also show that
\[
\mO{p, q}\cap \mU{p + q} =
\left\{
\left.
\amsmata{X_1}{0}{0}{X_2}\> 
\right|
\>
X_1\in\mO{p},\>\> X_2\in\mO{q}
\right\} \cong
\mO{p}\times \mO{q}.
\]
Therefore, we obtain the structure of $\mO{p, q}$:

\begin{prop}
\label{mOpqstruc}
The topological space $\mO{p, q}$ is homeomorphic
to $\mO{p}\times\mO{q}\times\reals^{pq}$.
\end{prop}

\medskip
Since $\mO{p}$ has two connected components when $p \geq 1$, we
see that $\mO{p, q}$ has four connected components when
$p, q\geq 1$. It is also obvious that
\[
\mSO{p, q}\cap \mU{p + q} =
\left\{
\left.
\amsmata{X_1}{0}{0}{X_2}\> 
\right|
\>
X_1\in\mO{p},\>\> X_2\in\mO{q},\>\> \det(X_1)\det(X_2) = 1
\right\}.
\]
This is a subgroup of $\mO{p}\times\mO{q}$ that we denote
$S(\mO{p}\times\mO{q})$. Furthermore, it is easy to show that
$\mfrac{so}(p, q) = \mfrac{o}(p, q)$.
Thus, we also have

\begin{prop}
\label{mSOpqstruc}
The topological space $\mSO{p, q}$ is homeomorphic
to $S(\mO{p}\times\mO{q})\times\reals^{pq}$.
\end{prop}

\medskip
Note that $\mSO{p, q}$ has two connected components
when $p, q\geq 1$. The connected component of
$I_{p+q}$ is a group denoted ${\bf SO}_0(p, q)$.
This latter space is homeomorphic to
$\mSO{p}\times\mSO{q}\times \reals^{pq}$.

\medskip
As a closing remark
observe that the dimension of all these spaces
depends only on $p + q$: It is $(p + q)(p + q - 1)/2$.

\bigskip\noindent
{\bf Acknowledgments}.
I thank Eric King whose incisive questions and
relentless quest for the ``essence'' of rotations
eventually caused a level of discomfort
high enough to force me
to improve the clarity of these notes.
Rotations are elusive! I also thank: Fred Lunnon 
for many insighful comments and for catching a number of typos;
and Troy Woo who also  reported some typos.